\documentclass[3p,times,12pt]{elsarticle}
\usepackage{graphicx}
\usepackage{longtable}
\usepackage{amsmath,amssymb,amsthm,amsfonts}
\usepackage{xcolor}

 \usepackage{epstopdf}

\usepackage[usenames,dvipsnames]{pstricks}

\usepackage{tikz}

\usetikzlibrary {positioning}
\definecolor {processblue}{cmyk}{0.96,0,0,0}

\newcommand{\Spec}{\mathrm{Spec}}

\newtheorem{theorem}{Theorem}[section]
\newtheorem{prop}[theorem]{Proposition}
\newtheorem{lemma}[theorem]{Lemma}
\newtheorem{exam}[theorem]{Example}
\newtheorem{conj}[theorem]{Conjecture}

\newtheorem{cor}[theorem]{Corollary}

\newcommand{\one}{\hbox{\it 1\hskip -3pt I}}
\newenvironment{pf}{\medskip {Proof:  \hspace*{-.4cm}}
    \enspace}{\hfill \qed \newline \medskip}

\begin{document}
\baselineskip 16pt
 \newcommand{\la}{\lambda}
 \newcommand{\si}{\sigma}
 \newcommand{\ol}{1-\lambda}
 \newcommand{\be}{\begin{equation}}
 \newcommand{\ee}{\end{equation}}
 \newcommand{\bea}{\begin{eqnarray}}
 \newcommand{\eea}{\end{eqnarray}}

\begin{frontmatter}
\title{Spectral Applications of Vertex-Clique Incidence Matrices Associated with a Graph}

\author{Shaun Fallat and  Seyed Ahmad Mojallal}
\address{Department of Mathematics and Statistics, University of Regina, Regina, Saskatchewan, S4S 0A2,  Canada}
\address{shaun.fallat@uregina.ca,\,ahmad\_mojalal@yahoo.com,}
\begin{abstract} 
In this paper, we demonstrate a useful interaction between the theory of clique partitions, edge clique covers of a graph, and the spectra of graphs. Using a clique partition and an edge clique cover of a graph we introduce the notion of a vertex-clique incidence matrix for a graph and produce new lower bounds for the negative eigenvalues and negative inertia of a graph. Moreover, utilizing these vertex-clique incidence matrices, we generalize several notions such as the signless Laplacian matrix, and develop bounds on the incidence energy and the signless Laplacian energy of the graph. 
More generally, we also consider the set $S(G)$ of all real-valued symmetric matrices whose off-diagonal entries are nonzero precisely when the corresponding vertices of the graph are adjacent. An important parameter in this setting is $q(G)$, and is defined to be the minimum number of distinct eigenvalues over all matrices in $S(G)$. For a given graph $G$ the concept of a vertex-clique incidence matrix associated with an edge clique cover is applied to establish several classes of graphs with $q(G)=2$.
\end{abstract}

\end{frontmatter}

 \section{Introduction}
 Let $G=(V,E)$ be a simple undirected graph with $n$ vertices and $m$ edges. A {\it clique} in $G$ is a subset $C \subseteq V$ such that all vertices in $C$ are adjacent.
An {\em edge clique cover} $F$ of  $G$ is a set of cliques $F=\{C_1, C_2, \ldots, C_k\}$ that together contain each edge of $G$ at least once.  The smallest size of an edge clique cover of $G$ is called the edge clique cover number of $G$ and is denoted by $cc(G)$. 
An edge clique cover of $G$ with size $cc(G)$ is referred to as a  {\it minimum edge clique cover} of $G$. A special case of an edge clique cover in which every edge belongs to exactly one clique is called a {\it clique partition} of $G$. The size of the smallest clique partition
 of $G$ is called the {\it clique partition number} of $G$, and is denoted by $cp(G)$. A clique partition of $G$ with size $cp(G)$ is referred to as a {\it minimum clique partition} of $G$. It is clear that both $cc(G)$ and $cp(G)$ exist as $E$ forms a clique partition (and hence an edge clique cover) of $G$. Further note that any minimum clique
 partition does not contain any cliques of size one, and, by convention, the clique partition number of the empty graph is defined to be zero. Information concerning clique partitions and edge clique covers of a graph can be found in the works \cite{Cav, Erd, LC, MCR}.

 Before defining the various matrices associated with a graph, we make note of the standard matrix notations: $I_n$ to denote the $n \times n$ identity matrix; $O$ to denote the zero matrix (size determined by context);  $J$ to denote the all ones matrix (size determined by context); and $\one$ to denote the all ones vector (size determined by context).


 Given a graph $G$ with $V=\{1,2, \ldots, n\}$ and $E=\{e_1, e_2, \ldots, e_m\}$, the {\it (vertex-edge) incidence matrix} $M$ of $G$ is the $n\times m$ matrix defined as follows: the rows and the columns of $M$ are indexed by $V$ and $E$, respectively; and the $(i,j)$-entry of $M$ is 0 if $i \not\in e_j$ and 1 otherwise. Similarly, the
 adjacency matrix ${\cal A}={\cal A}(G)=(a_{ij})$ is a  $(0, 1)$-matrix of $G$ such that $a_{ij}=1$ if $ij\in E(G)$ and 0 otherwise. It is well-known that \cite{GR}
 \begin{equation}\label{equ5}
 M M^T=Q(G),~~~\mbox{and}~~~M^T M={\cal A}(L_G)+2I_m,
 \end{equation}
 where
 $D(G)$ is the diagonal matrix of vertex degrees ($d_i=\deg(i)$, $i=1,2,\ldots, n$) and the matrix $Q(G)= D(G)+{\cal A}(G)$
 is known as the {\it signless Laplacian matrix} of the graph $G$; the {\it line graph}, $L_G$, of the graph $G$ is the graph whose vertex set is in one-to-one correspondence with the set of edges of
 $G$, where two vertices of $L_G$ are adjacent if and only if the corresponding edges in $G$ have a vertex in common \cite{Har}.
Finally, the equations in  (\ref{equ5}) imply an important spectral relation between the signless Laplacian matrix $Q(G)$ and ${\cal A}(L_G)$, see Lemma \ref{line}.

 As we are also interested in studying more general symmetric matrices associated to a graph on $n$ vertices, we let $S(G)$ denote the collection of real symmetric matrices  $A=(a_{ij})$ such that for $i\neq j$, $a_{ij}\neq 0$ if and only if $ij\in E(G)$. The main diagonal entries of any such $A$ in $S(G)$ are not constrained. Observe that for any graph $G$, both $Q(G)$ and ${\cal A}(G)$ belong to $S(G)$.
 
 We denote the spectrum of $A$, i.e., the multiset of eigenvalues of $A$, by $\Spec(A)$. In particular, $$\Spec(A)=\{\lambda_1^{[m_1]},\, \lambda_2^{[m_2]},\, \ldots,\,\lambda_q^{[m_q]}\},$$ where
 the distinct eigenvalues of $A$ are given by $\lambda_1< \lambda_2< \cdots <\lambda_q$ with corresponding multiplicities of these eigenvalues are $m_1, m_2, \ldots, m_q$ respectively. Further we consider the ordered multiplicity list of $A$ as the sequence $m(A) =(m_1, m_2, \ldots, m_q)$. For brevity, a simple eigenvalue $\lambda_k^{[1]}$ is simply denoted by  $\lambda_k$.
 
 Given a graph $G$, the spectral invariant $q(G)$ is defined as follows:
 $$q(G)=\min\{q(A)\,:\, A\in S(G)\},$$ where $q(A)$ is the number of distinct eigenvalues of $A$ (see \cite{AACF, LOS}).
The spectral invariant $q(G)$ is called the {\em minimum number of distinct eigenvalues of the graph $G$}. The class of matrices $S(G)$ has been of interest to many researchers recently (see \cite{Fal, FH, Fal, Fer} and the references therein), and there has been considerable development on the inverse eigenvalue problem for graphs (see \cite{Hog}) which continues to receive considerable and deserved attention, as it remains one of the most interesting unresolved issues in combinatorial
matrix theory. Recently, J.  Ahn et al. \cite{AAB} offered a complete solution to the ordered multiplicity inverse eigenvalue problem for graphs on six vertices. 

Using the notions of clique partitions and edge clique covers of a graph we generalize  the conventional vertex-edge incidence matrix $M$ by considering a new incidence matrix called the {\it vertex-clique incidence matrix} of a graph. Suppose $F=\{C_1, C_2, \ldots, C_k\}$ is an edge clique cover of a graph $G$ with $V=\{1,2,\ldots, n\}$. The vertex-clique incidence matrix $M_F$ of $G$ associated with the edge clique cover $F$ is defined as follows:
the $(i,j)$-entry of $M_F$ is real and nonzero if and only if the vertex $i$ belongs to the clique $C_j\in F$. In the particular case when $F$ is actually a clique partition, the vertex-clique incidence matrix, in this case, is denoted by ${\cal M}_F$, and the $(i,j)$-entry of ${\cal M}_F$ is equal to one if and only if the vertex $i$ belongs to the clique $C_j\in
 F$. We observe that for any graph $G$ the vertex-clique incidence matrix corresponding to a clique partition $F$, preserves several main properties of its vertex-edge incidence matrix.
 For instance, in Section \ref{section:matrices}, ${\cal M}_F {\cal M}^T_F={\cal D}_F+{\cal A}$, where ${\cal D}_F=diag(t_1^F, t_2^F, \ldots, t_n^F)$ with $t_i^F\le d_i$, where both sequences $t_i^F$ and $d_i$ are in non-increasing order.
 This fact enables us to determine new lower bounds for the negative eigenvalues of the graph.


 The paper is organized as follows. In Section~\ref{section:notation}, we provide the necessary notions, notations, and known results that are needed in the sections containing our main observations. In Section~\ref{section:matrices}, using the notion of a clique partition $F$ of a graph $G$, we define signless Laplacian matrix of the graph $G$ associated with the clique partition $F$. A new graph $P_G$ is introduced as a generalization for the line
 graph of $G$. In Subsection~\ref{section:spectrum}, applying this new theory of a vertex-clique incidence matrix, we produce lower bounds for the negative eigenvalues of the graph. Moreover, we present lower bounds for the negative inertia
 $\nu^-(G)$ of a graph $G$ in terms of its order $n$ and the rank of its vertex-clique incidence matrix. We also provide a sufficient condition under which the well-known inequality $\nu^-(G)\le n-\alpha(G)$ holds with equality, where $\alpha(G)$ is the independence number
 of $G$. In Subsection~\ref{section:energy}, we introduce new graph energies associated with a clique partition $F$ of the graph $G$ and study several associated properties. Moreover, new upper bounds for the energies of the graph $G$ and its clique partition graph and line graph are determined. In Section~\ref{Vertex-clique incidence matrix}, studies on the vertex-clique incidence matrix of a graph associated with an edge clique cover lead to a derivation of some new classes of graphs with $q(G)=2$ (see also Subsection~\ref{q=2}).

 \section{Notations and preliminaries}\label{section:notation}

 In this section, we list some known notions, notations, and results that are needed in the remaining
 sections. 
 
 We start this section by introducing the notion of the eigenvalues of a graph. The eigenvalues $\lambda_1,\,\lambda_2,\ldots,\,\lambda_n$ of the adjacency matrix ${\cal  A}(G)$ (or shortened to ${\cal A}$ when reference to the graph $G$ is clear from context) of the graph $G$ are also
 called the {\it eigenvalues of $G$}.  The number of positive (negative) eigenvalues in the spectrum of the graph $G$ is called the {\it positive (negative) inertia} of the graph $G$, and is denoted by $\nu^+(G)$ ($\nu^-(G)$).
 The {\it energy} of the graph $G$ is defined as
 \begin{equation}
 {\cal E}(G)=\sum\limits^n_{i=1}\,|\lambda_i|\,.\label{bm0}
 \end{equation}
 Further details on various properties of graph energy can be found in \cite{GU1,GU2,KM,LSG,MC}.
 Suppose $q_1,\,q_2,\ldots,\,q_n$ be the eigenvalues of the matrix $Q(G)$.
 Then the {\it signless Laplacian energy} of the graph $G$ is defined as \cite{Ab}
 \begin{equation}                    \label{bm2}
 LE^+ = LE^+(G)=\sum^n_{i=1}\big|q_i-\frac{2m}{n}\big| .
 \end{equation}
More information on properties of the signless Laplacian energy can be found in \cite{Ab}, and 
 the energy of a line graph and its relations with other graph energies are studied in \cite{AH5,GRMC}.

A {\em subgraph} $H$ of a graph $G$ is a graph whose vertex set and edge set are subsets of those of $G$. If $H$ is a subgraph of $G$, then $G$ is said to be a {\em supergraph} of $H$. The subgraph of $G$ obtained by deleting either a vertex $v$ of $G$ or an edge $e$ of $G$ is denoted by $G-v$ and $G-e$, respectively. Suppose $H$ is a graph on $n$ vertices. Then we let $K_n \backslash H$ denote the graph obtained from the complete graph, $K_n$, by removing the edges from $H$. An {\it independent set} in the graph $G$ is a set of vertices in $G$, no two of which are adjacent. The {\it independence number}  $\alpha(G)$ of $G$ is the number of vertices in a largest independent set of $G$. A {\em matching} in a graph $G$, is simply a collection of independent edges from $G$ (i.e., no two edges in a  matching share a common vertex from $G$). Additionally, a matching is referred to as {\em perfect} if each vertex from $G$ is incident with one edge from the matching.

 An $n \times n$ symmetric real matrix $B$ is a positive semi-definite matrix if all of its eigenvalues are nonnegative. In this case, we denote $B\ge 0$. For real symmetric matrices $B$ and $C$, if $B-C\ge 0$, then we write $B\ge C$.
 \begin{lemma}{\rm  \cite{Ber} } \label{lem3}
 Let $A$ and $B$ be Hermitian matrices of order $n$, and assume that $A\le B$. Then for all $i=1, 2, \ldots, n$,
 $$\lambda_i(A)\le \lambda_i(B),$$ where $\lambda_i(M)$ is the $i$th largest eigenvalue of a square matrix $M$.
 \end{lemma}

 The following result was obtained in \cite{GR}.
 \begin{lemma}\label{lem1}{\rm \cite{GR}}
 If $B$ and $C$ are matrices such that $BC$ and $CB$ are both defined, then $BC$ and $CB$ have the same nonzero eigenvalues with the same multiplicity.
 \end{lemma}

Let $\circ$ denote the Schur (also known as the Hadamard or entry-wise) product. 
The $n \times n$ symmetric matrix $A$ has the {\it Strong Spectral Property} (or $A$ has the SSP for short) if the only symmetric matrix $X$ satisfying $A\circ X=O$, $I\circ X=O$ and $[A,\, X]=AX-XA=O$ is $X=O$ (see \cite{BFH}).
The following result is given in \cite[Thm. 10]{BFH}.
\begin{lemma}{\rm \cite{BFH}}\label{ssp}
If $A\in S(G)$ has the SSP, then every supergraph of $G$ with the same vertex set has a matrix realization that has the same spectrum as $A$ and has the SSP.
\end{lemma}

Given two graphs $G$ and $H$, the join of $G$ and $H$, denoted by $G \vee H$, is the graph obtained from $G \cup H$, by adding all possible edges between $G$ and $H$.
Suppose $G$ is a graph with $q(G)=2$. Then among all matrix realizations $A$ in $S(G)$ with two distinct eigenvalues, we define the multiplicity bi-partition $[n-k, \, k]$ associated to $A$ if the two eigenvalues of $A$ has respective multiplicities $n-k$ and $k$. Further we  define the minimal multiplicity bi-partition $MB(G)$ to be the least integer $k\le \lfloor \frac{n}{2} \rfloor$ such that $G$ achieves
the multiplicity bi-partition $[n-k, \, k]$.
We close this section with two useful results concerning specific classes of graphs realizing two distinct eigenvalues with respect to the set $S(G)$. 

 \begin{lemma}{\rm \cite{BHL,CGMJ} }\label{MB}
Let $G$ be a connected graph on $n$ vertices. Then\\
(1) $MB(G) = 1$ if and only if $G$ is the complete graph, $K_n$.\\
(2) $MB(G) = 2$ if and only if $$G=(K_{p_1}\cup K_{q_1})\vee (K_{p_2}\cup K_{q_2})\vee \ldots (K_{p_k}\cup K_{q_k})$$ for non-negative integers $p_1, \ldots, p_k,q_1, \ldots, q_k$ with $k > 1$, and $G$ is not isomorphic to either one of a complete graph or $G=(K_{p_1}\cup K_{q_1})\vee K_1$. 
\end{lemma}

 \begin{lemma}{\rm \cite{LOS2} }\label{lemma-join}
If $G$ is a connected graph of order $n \in \{l, l+ 1, l+ 2\}$ and 
$n_1,  \ldots, n_l \in \mathbb{N}$, then $q(G\vee \cup_{j\in [l]} K_{n_j} ) = 2$.
\end{lemma}

 \section{Matrices associated with a clique partition}\label{section:matrices}
 
  In this section, we use of the vertex-clique incidence matrix associated with a clique partition of a graph $G$.
  Recall that for a graph $G=(V, E)$ with the vertex set $V=[n]=\{1,2, \ldots, n\}$ and $m=|E|$ edges, and for a given clique partition $F=\{C_1, C_2, \ldots, C_k\}$ of $G$, consider the matrix ${\cal M}_F$ with rows and columns indexed by the vertices in $V$ and the cliques in $F$, respectively, such that
 the $(i,j)$-entry of ${\cal M}_F$ is equal to one if and only if the vertex $i$ belongs to the clique $C_j\in
 F$. Observe that when $F=E$, ${\cal M}_F$ as simply the conventional incidence matrix of the graph $G$. 
 For each vertex $i\in [n]$ of the graph $G$, we define a new parameter $t_i^F=t_i^F(G)$ to be the number of cliques in $F$ containing the vertex $i$, that is,
 $$t_i^F=|\{j\in [k]\,:\, C_j\in F,\,i\in C_j\}|.$$ We call $t_i^F(G)$ {\it the clique-degree} of the vertex $i$ in graph $G$ associated with
 $F$, and, without loss of generality, we assume that $t_1^F\ge t_2^F \ge \ldots \ge t_n^F$.
 Given clique partition $F=\{C_1, C_2, \ldots, C_k\}$ of $G$, we consider different possible classes of graphs as
 follows:\\
 $(i)$ The graph $G$ is {\it $t$ clique-regular} if $t_1^F=\cdots=t_n^F=t$,\\
 $(ii)$  The graph $G$ is {\it $s$ clique-uniform} if $|C_1|=\cdots=|C_k|=s$,\\
 $(iii)$ The graph $G$ is  {\it $(s,t)$ regular} if $t_1^F=\cdots=t_n^F=t$ and
 $|C_1|=\cdots=|C_k|=s$.

 Any graph is 2 clique-uniform and any $d$-regular graph is also $d$ clique-regular using the trivial clique partition $F=E$.

 Let ${\cal D}_F$ be the $n \times n$ diagonal matrix with row and column indexed by the vertex set
 $V$ with $(i,i)$-entry equal to $t_i^F$, that is, ${\cal D}_F=diag(t_1^F,\ldots, t_n^F)$. The inner product of any two distinct rows of ${\cal M}_F$ indexed by vertices $i$ and $j$
 is equal to the number of cliques in $F$ containing the vertices $i$ and $j$. By definition of the clique partition $F$, if $i$ and $j$ are adjacent, then this number is equal to 1 and otherwise 0.
 This leads to the following result:
 \begin{theorem}
 Let ${\cal M}_F$ be the vertex-clique incidence matrix of $G$ associated with a given clique partition
 $F$. Then ${\cal M}_F{\cal M}_F^T={\cal D}_F+{\cal A},$  where ${\cal D}_F=diag(t_1^F,\ldots, t_n^F)$ and ${\cal A}$ is the adjacency matrix of $G$.
 \end{theorem}

 As mentioned above, in the case of $F=E$, the matrix ${\cal M}_F$ is the incidence matrix $M$ of $G$ and consequently, ${\cal M}_F{\cal M}^T_F=M M^T$
 is the signless Laplacian matrix of $G$, where we assume that the sequence of vertex degrees is ordered as $d_1 \ge d_2 \ge \cdots \ge d_n$. Notice that in this case, $t_i^F=d_i$ for $1\le i \le n$. Motivated by this observation, for any clique partition $F$ we call ${\cal Q}_F={\cal M}_F{\cal M}_F^T$ the {\it signless Laplacian matrix of the graph $G$ associated with the clique partition $F$}.
 Since we always have $D\ge {\cal  D}_F$, it follows $Q=D+{\cal A}\ge {\cal  D}_F+{\cal A}={\cal Q}_F\ge 0$.  Now define {\it the clique partition graph} $P_G$ with $k$ vertices, where each vertex $i$ corresponds to each clique $C_i$ in $F$ such that each pair of vertices of $P_G$ are adjacent if and only if the corresponding cliques in $F$ have a vertex in
  common. If $F=E$, then $P_G=L_G$, that is, the line graph of $G$.
  The inner product of two columns of ${\cal M}_F$ is nonzero if and only if the corresponding cliques have a common vertex. From the definition of a clique partition, this nonzero value must be 1. These facts immediately yield the following result:
 \begin{theorem}
 Let ${\cal M}_F$ be the incidence matrix of $G$ associated with a clique partition
 $F$. Then ${\cal M}_F^T {\cal M}_F={\cal S}_F+{\cal A}(P_G),$  where ${\cal S}_F=diag(s_1^F,\ldots, s_k^F)$ and $s_i^F=|C_i|$ and ${\cal A}(P_G)$ stands for the adjacency matrix of the graph $P_G$.
 \end{theorem}

 For the case of $F=E$, we have ${\cal M}^T_F{\cal M}_F=M^T M=2I_m+{\cal A}(L_G)$, and $P_G=L_G$ so $s_i^F=2$ for $1\le i \le k=m$.

\subsection{Applications of the vertex-clique incidence matrix to graph spectrum}\label{section:spectrum}

 In this section, we develop several results on the spectrum of the graph
 $G$ and its clique partition graph $P_G$ by the vertex-clique incidence matrix of a graph.
Considering ${\cal R}_F={\cal M}^T_F{\cal M}_F$ with Lemma \ref{lem1} we conclude that the nonzero eigenvalues of matrices ${\cal Q}_F$ and ${\cal R}_F$ are the same. This fact leads to the following basic results.
 \begin{theorem}\label{pro2} We have the following. \\
 $(i)$\, If $1\le i \le \min\{n, k\}$, then $\lambda_i({\cal Q}_F)=\lambda_i({\cal
 R}_F)$.\\
 $(ii)$\, If $\min\{n, k\}=n$ then  $\lambda_i({\cal R}_F)=0$ for $n+1\le i\le k$.\\
 $(iii)$\, If $\min\{n, k\}=k$ then $\lambda_i({\cal Q}_F)=0$ for $k+1\le i\le n$.
 \end{theorem}

 Recall that if $F=E$, then ${\cal
 Q}_F=Q$ and ${\cal R}_F=2I_m+{\cal A}(L_G)$. Combining these equations with Theorem  \ref{pro2} leads to the following well-known
 result \cite{cv, AH5}:
 \begin{lemma}\label{line} Let $G$ be a graph of order $n$ with $m$ edges. Then
  $$q_i(G)=2+\lambda_i(L_G)~~~\mbox{for}~1\le i \le \min\{n,m\}.$$
  In particular if $m>n$ then $\lambda_i(L_G)=-2$ for $i>n$, and if $n>m$ then   $q_i(G)=0$ for $i> m$.
 \end{lemma}
 The following result is obtained by applying Theorem \ref{pro2} for a $(s, t)$ regular graph $G$ with the clique partition $F$.
 \begin{theorem}\label{pro1}
 Let $G$ be a $(s, t)$ regular graph of order $n$ with a clique partition $F$ of size $k$.\\
 $(i)$\, If $1\le i \le \min\{n, k\}$ then
 $\lambda_i(G)-\lambda_i(P_G)=s-t.$\\
 $(ii)$\, If $\min\{n, k\}=n$ then $\lambda_i(P_G)=-s$ for $n+1\le i\le
 k$.\\
 $(iii)$\, If $\min\{n, k\}=k$ then $\lambda_i(G)=-t$ for $k+1\le i\le n$.
 \end{theorem}
 \begin{pf}
 $(i)$\, By Theorem \ref{pro2} $(i)$, if $1\le i \le \min\{n, k\}$, then $\lambda_i({\cal Q}_F)=\lambda_i({\cal
 R}_F)$, that is, $\lambda_i({\cal D}_F+{\cal A}(G))=\lambda_i({\cal S}_F+{\cal A}(P_G))$, that is, $\lambda_i(tI_n+{\cal A}(G))=\lambda_i(sI_k+{\cal A}(P_G))$, that is, $t+\lambda_i(G)=s+\lambda_i(P_G).$

\noindent
 $(ii)$\, By Theorem \ref{pro2} $(ii)$, if $\min\{n, k\}=n$ then  $\lambda_i({\cal R}_F)=0$ for $n+1\le i\le k$, that is, $\lambda_i(sI_k+{\cal A}(P_G))=0$ for $n+1\le i\le k$, that is, $\lambda_i(P_G)=-s$ for $n+1\le i\le k$.

\noindent
 $(iii)$\, By Theorem \ref{pro2} $(iii)$, if $\min\{n, k\}=k$ then $\lambda_i({\cal Q}_F)=0$ for $k+1\le i\le n$, that is, $\lambda_i(tI_n+{\cal A}(G))=0$ for $k+1\le i\le n$, that is, $\lambda_i(G)=-t$ for $k+1\le i\le n$.
 \end{pf}

 \begin{exam} {\rm $(i)$\, Considering the complete graph $K_n$ and its minimum clique partition $F$ with only one clique, we have ${\cal M}_F=\one_n$, ${\cal M}_F {\cal M}_F^T=J_n$ and ${\cal M}_F^T {\cal M}_F=[n]$.
 Applying Theorem \ref{pro1} here we have
 $t_i^F=1$ for $1\le i \le n$, $k=1$ and $s_1^F=n$, that is, $K_n$ is a $(n, 1)$ regular graph. From this with Theorem \ref{pro1} $(i)$ we arrive at $1+\lambda_1(K_n)=n+\lambda_1(K_1),$ that is, $\lambda_1(K_n)=n-1$,
 and by Theorem \ref{pro1} $(iii)$, $\lambda_i(K_n)=-1$ for $2\le i \le n$.\\
 $(ii)$\, Considering the clique
 partition  $$F=\Big\{C_1=\{1,2,6\}, C_2=\{2, 3, 4\}, C_3=\{1,3,5\},
 C_4=\{4,5,6\}\Big\}$$ for $G$ isomorphic to the complete tripartite graph $K_{2,2,2}$  (or $G\cong K_{2,2,2}$) in Figure \ref{k222}, we have
 $s_i^F=3$ for $i\in [4]$ and $t_j^F=2$ for $j\in [6]$. Then $G$ is a $(3, 2)$ regular graph.
 Moreover, $${\cal M}_F=\left(
                           \begin{array}{cccc}
                             1 & 0 & 1 & 0 \\
                             1 & 1 & 0 & 0 \\
                             0 & 1 & 1 & 0 \\
                             0 & 1 & 0 & 1 \\
                             0 & 0 & 1 & 1 \\
                             1 & 0 & 0 & 1 \\
                           \end{array}
                         \right),~{\cal Q}_F=\left(
                                      \begin{array}{cccccc}
                                        2 & 1 & 1 & 0 & 1 & 1 \\
                                        1 & 2 & 1 & 1 & 0 & 1 \\
                                        1 & 1 & 2 & 1 & 1 & 0 \\
                                        0 & 1 & 1 & 2 & 1 & 1 \\
                                        1 & 0 & 1 & 1 & 2 & 1 \\
                                        1 & 1 & 0 & 1 & 1 & 2 \\
                                      \end{array}
                                    \right),~{\cal R}_F=\left(
                                                            \begin{array}{cccc}
                                                              3 & 1 & 1 & 1 \\
                                                              1 & 3 & 1 & 1 \\
                                                              1 & 1 & 3 & 1 \\
                                                              1 & 1 & 1 & 3 \\
                                                            \end{array}
                                                          \right)$$ and by
 Theorem \ref{pro1}, we have
 $\lambda_i(G)=1+\lambda_i(P_G)$ for $1\le i \le 4$ and $\lambda_i(G)=-2$ for $i=5,6$. From these facts with $P_G\cong K_4$, we
 arrive at $\Spec(G)=\{4, 0, 0, 0, -2, -2\}$.
 \begin{figure}[htb] 
 \begin{center}
 \begin{tikzpicture}
  [scale=.6,auto=left,every node/.style={circle,fill=black!15}]
  \node (n6) at (0,10) {6};
  \node (n4) at (8,10) {4};
  \node (n5) at (4,3)  {5};
  \node (n1) at (3,7)  {1};
  \node (n2) at (4,9)  {2};
  \node (n3) at (5,7)  {3};

 \foreach \from/\to in {n6/n4,n6/n5,n6/n2,n6/n1,n1/n2,n1/n3,n2/n3,n5/n1,n5/n3,n5/n4,n4/n3,n4/n2}
 \draw (\from) -- (\to);
 \end{tikzpicture}
 \caption{The graph $G\cong K_{2,2,2}$.}
 \label{k222}
 \end{center}
 \end{figure}
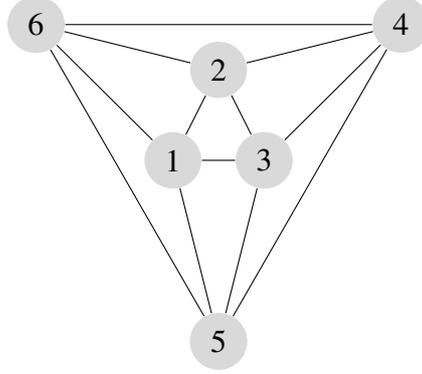}
 \end{exam}

 Now applying theory of clique partitions and vertex-clique incidence matrices, we obtain a new lower bound for the smallest eigenvalue of a graph.
 \begin{theorem}\label{thm1}
 Let $G$ be a graph of order $n$ and let $t_1^F$ be the largest clique-degree of $G$ with a given clique partition $F$. Then
 \begin{equation}\label{eq2}
 \lambda_n(G)\ge -t_1^F.
 \end{equation}  Moreover, if equality holds in (\ref{eq2}), then
 $rank({\cal M}_F)<n$ and if $rank({\cal M}_F)<n$ and $G$ is clique-regular, then equality holds in (\ref{eq2}).
 \end{theorem}
 \begin{pf}
 Since ${\cal Q}_F={\cal D}_F+{\cal A}$ is a positive semi-definite matrix, we have ${\cal D}_F\ge -{\cal A}$ and by Lemma \ref{lem3} we arrive at
 \begin{equation}\label{eq4}
 \lambda_i({\cal D}_F)\ge \lambda_i(-{\cal A})~~\mbox{for}~1\le i \le n.
 \end{equation} Considering $i=1$ we arrive at $-\lambda_n(G)=\lambda_1(-{\cal A})\le \lambda_1({\cal D}_F)=t_1^F$, which gives the required result in (\ref{eq2}).

 For the second part of the proof, suppose that $\lambda_n(G)=-t_1^F$. Then $\lambda_n(t_1^FI+{\cal A})=0$. This with the relation $0\le {\cal Q}_F={\cal D}_F+{\cal A}\le t_1^FI+{\cal A}$, gives $\lambda_n({\cal Q}_F)=0$, that is, $rank({\cal M}_F)=rank({\cal Q}_F)<n$.
Now we assume that  $t_1^F=\cdots=t_n^F$. If $rank({\cal M}_F)<n$ then $rank({\cal Q}_F)<n$, that is, $\lambda_i({\cal Q}_F)=0$ for $1+k\le i\le n$, that is, $t_1^F+\lambda_i(G)=0$ as ${\cal Q}_F=t_1^F I+{\cal A}$, that is, $\lambda_n(G)=-t_1^F$
 with the multiplicity at least $n-k$.
 \end{pf}
 \begin{cor}
 All regular bipartite graphs and all clique-regular graphs with $n> |F|$ satisfy the equality in (\ref{eq2}).
 \end{cor}
 \begin{pf} First we assume that $G$ is a regular bipartite graph.
 Since $G$ is bipartite, we have $t_i^F=d_i$ for $i\in [n]$ and
 $q_n=\lambda_n(Q)=0$. On the other hand, since $G$
 is regular, we have $t_1^F=\cdots=t_n^F$. These facts with Theorem
 \ref{thm1} gives the fact that all regular bipartite graphs
 satisfies the equality in (\ref{eq2}).

 Next assume that $G$ is a clique-regular graph with $n>k=|F|$. Since $rank({\cal M}_F)\le \min{\{n,k\}}\le k<n$, the desired result is obtained by Theorem
 \ref{thm1}.
 \end{pf}

 Theorem \ref{thm1} holds for any clique partition $F$ of $G$, which leads to the following
 result:
 \begin{cor}
 Let $G$ be a graph of order $n$ and let $t_1^F$ be the largest clique-degree of $G$ with a given clique partition $F$. Then
 \begin{equation}\nonumber
 \lambda_n(G)\ge -\min_{F}t_1^F,
 \end{equation}  where the minimum is over all clique partitions $F$ of $G$.
 \end{cor}

 The following example shows that for the equality $\lambda_n(G)=-t_1^F$
 the graph $G$ does not need to be clique-regular.
 \begin{exam}
 For the graph $G$ given in Figure \ref{t1eq}, we have $$F=\Big\{\{1,2\}, \{2,3\}, \{1,3,4,6,7\}, \{4,5\}, \{5,6\}\Big\}.$$
 This gives $t_i^F=2$ for $i\in [6]$ and $t_7^F=1$. The graph is the line graph of the  graph $H\cong K_1\vee (2K_2\cup K_1)$ of order $6$ with $7$ edges. Then the smallest eigenvalue of $G$ is $\lambda_7(G)=\lambda_7(L_H)=-2=-t_1^F$
 while $t_1^F\neq t_7^F$.
  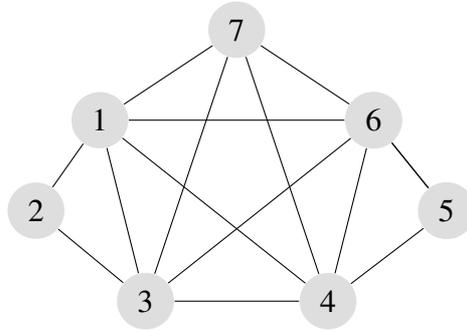
\begin{figure}[htb]
  \begin{center}
 \begin{tikzpicture}
 [scale=.6,auto=left, every node/.style={circle,fill=black!13}]
  \node (n6) at (-.4,5) {2};
  \node (n4) at (8.6,5) {5};
  \node (n5) at (2,3)  {3};
  \node (n1) at (1,7)  {1};
  \node (n2) at (4,9)  {7};
  \node (n3) at (7,7)  {6};
  \node (n7) at (6,3)  {4};

  \foreach \from/\to in {n5/n2,n6/n5,n6/n1,n1/n2,n1/n3,n2/n3,n2/n7,n5/n1,n5/n3,n4/n3,n4/n3,n7/n1,n7/n3,n7/n5,n7/n4}
    \draw (\from) -- (\to);
 \end{tikzpicture}
 \caption{The Graph $G$.}
 \label{t1eq}
 \end{center}
\end{figure}
 \end{exam}

 In the following we provide a lower bound for the negative inertia of
 a graph $G$ of order $n$.
 \begin{theorem}\label{thm3}
 Let $G$ be a graph of order $n$. Then
 \begin{equation}\label{equ1}
 \nu^-(G)\ge n-\min_{F}{rank({\cal M}_F)},
 \end{equation} where minimum is over all clique partitions $F$ of $G$. Moreover, if $\min_{F}{rank({\cal M}_F)}<n$, then $-t_1^F\le \lambda_i(G)\le -t_n^F$ for $1+\min\limits_{F}{rank({\cal M}_F)}\le i \le n$.
 \end{theorem}
 \begin{pf}
 If $\min\limits_{F}{rank({\cal M}_F)}=n$, then the result in (\ref{equ1}) is obvious. Assume that $F_1$ is a clique partition of $G$ with $rank({\cal M}_{F_1})=\min\limits_{F}{rank({\cal M}_F)}<n$.
 In this case, since $rank({\cal Q}_{F_1})=rank({\cal M}_{F_1})$ and ${\cal Q}_{F_1}$ is positive semi-definite matrix, we have $\lambda_i({\cal Q}_{F_1})=0$ for $1+rank({\cal M}_{F_1})\le i \le n$.
 From this and the fact that $t_n^F+\lambda_i(G)\le \lambda_i({\cal Q}_{F_1})\le t_1^F+\lambda_i(G)$, we have $-t_1^F\le \lambda_i(G)\le -t_n^F<0$  for $1+rank({\cal M}_{F_1})\le i \le n$, which gives the desired results.
 \end{pf}

 The following result is obtained by Theorem \ref{thm3} and the fact $rank({\cal M}_F)\le |F|$.
 \begin{cor}\label{thm2}
 Let $G$ be a graph of the order $n$ and a clique partition $F$ such that
 $n>|F|$. Then \\
 $(i)$\, $-t_1^F\le \lambda_i(G)\le -t_n^F$ for $|F|+1\le i \le n$.\\
 $(ii)$\, $\nu^-(G)\ge n-|F|$.
 \end{cor}

 Considering $F$ as a minimum clique partition of $G$, we arrive at the following result:
 \begin{cor}\label{cor1}
 Let $G$ be a graph of the order $n$ and clique partition number $cp(G)$. If $cp(G)<n$, then \\
 $(i)$\, $-t_1^F\le \lambda_i(G)\le -t_n^F$ for $cp(G)+1\le i \le n$.\\
 $(ii)$\, $\nu^-(G)\ge n-cp(G)$.
 \end{cor}

 \vspace{3mm}

 For any graph $G$ of order $n$ we have \cite{cv}
 \begin{equation}\label{alpha}
 \alpha(G)\le \min\{n-
 \nu^-(G),\, n-\nu^+(G)\},
 \end{equation} where $\nu^-$ and $\nu^+$ are the negative and positive parts of the inertia, respectively of the graph
 $G$. This implies that \begin{equation}\label{inertia}\nu^-(G)\le n-\alpha(G). \end{equation} In the following we
 give a sufficient condition under which the equality in
 (\ref{inertia}) holds.
 \begin{theorem}  Let $G$ be a graph of order $n$ with the independence number $\alpha(G)$ and the clique partition number
 $cp(G)$. If $F$ is a clique partition with $rank({\cal M}_F)=\alpha(G)$, then $\nu^-(G)=n-\alpha(G).$ In particular, if $cp(G)=\alpha(G)$, then $\nu^-(G)=n-\alpha(G).$
 \end{theorem}
 \begin{pf}
 By Theorem \ref{thm3} we have $$\nu^-(G)\ge
 n-rank({\cal M}_F)=n-rank({\cal Q}_F)=\eta({\cal Q}_F).$$
 This fact along with (\ref{inertia}) gives
 \begin{equation}\label{equ2}
 \eta({\cal Q}_F)\le \nu^-(G)\le n-\alpha(G).
 \end{equation} The assumption that ${\rm rank({\cal M}_F)}=\alpha(G)$ is
 equivalent to $\eta({\cal Q}_F)=n-\alpha(G)$. This with (\ref{equ2})
 gives the first required result.

 Without loss of generality, we may assume that the vertex set $[\alpha]$ is a maximum independent set in $G$ and $C_i$ is a clique of a minimum clique partition $F_m$ containing the vertex $i\in [\alpha]$.
 Now in  ${\cal M}_{F_m}$ we consider the submatrix induced by the rows and columns corresponding to the vertex set $[\alpha]$ and the clique set $\{C_i\,:\, i\in
 [\alpha]\}$, respectively. Obviously, this square principal submatrix is equivalent to the identity matrix of size $\alpha$ and hence ${\rm rank}({\cal Q}_{F_m})\ge {\rm rank}(I_{\alpha})=\alpha$. Since
 ${\rm rank}({\cal Q}_{F_m})\le cp(G)$ and using the assumption $cp(G)=\alpha(G)$ we arrive at ${\rm rank}({\cal Q}_{F_m})=\alpha={\rm rank} ({\cal M}_{F_m})$ and therefore $\nu^-(G)=n-\alpha(G)$ by the first part of the theorem.
 \end{pf}

 The following result is obtained from (\ref{eq4}).
 \begin{theorem} \label{thm6} Let $G$ be a graph of order $n$ and the negative inertia $\nu^-$. Let $t_i^F$ be the $i$th largest clique-degree of $G$ with a clique partition $F$.  Then for $1\le i \le \nu^-$, we have
 \begin{equation}\label{eq3}
 \lambda_{n-i+1}(G)\ge -t_i^F.
 \end{equation} Equality holds in  (\ref{eq3}) if $G$ is a clique-regular graph with $\nu^-=n-|F|$.
 \end{theorem}

 Since ${\cal R}_F$ is a positive semi-definite matrix, by a similar manner used in the proof of Theorem \ref{thm1}, we obtain the following result.
 \begin{theorem}\label{thm10}
 Let $G$ be a graph of order $n$ with a clique partition $F=\{C_1, \ldots, C_k\}$ and let $|C_i|=s_i^F$ for $1\le i\le k$ such that
 $s_1^F\ge s_2^F \ge \ldots \ge s_k^F$. Then 
 \begin{equation}\label{equ80}\lambda_k(P_G)\ge -s_1^F.
 \end{equation}
 Equality holds in (\ref{equ80}) if  $G$ is a $s_1^F$ clique-uniform graph with $k>n$.
 \end{theorem}
 \begin{pf}
 Since ${\cal R}_F={\cal S}_F+{\cal A}(P_G)$ is a positive semi-definite matrix, we have ${\cal S}_F\ge -{\cal A}(P_G)$ and by Lemma \ref{lem3}, it follows that
 \begin{equation}\label{eq5}
 \lambda_i({\cal S}_F)\ge \lambda_i(-{\cal A}(P_G))~~\mbox{for}~1\le i \le k.
 \end{equation} Considering $i=1$ we have $-\lambda_k(P_G)=\lambda_1(-{\cal A}(P_G))\le \lambda_1({\cal
 S}_F)=s_1^F,$ which gives the required result in (\ref{equ80}).

 Now assume that $G$ is a $s_1^F$ clique-uniform graph with $k>n$. By Theorem \ref{pro2} (ii) with $k>n$, we arrive at $\lambda_i({\cal R}_F)=0$ for
 $n+1\le i \le k$. On the other hand, since $s_1^F=\cdots=s_k^F$ we have ${\cal R}_F=s_1^F I_k+{\cal A}(P_G)$, and consequently $\lambda_i({\cal R}_F)=s_1^F+\lambda_i(P_G)=0$.
 That is, $\lambda_i(P_G)=-s_1^F$ for $n+1\le i \le k$, that is, $\lambda_k(P_G)=-s_1^F$ with  multiplicity at least $k-n$.
 \end{pf}

 Theorem \ref{thm10} holds for any clique partition $F$ of $G$, which gives the following
 result:
 \begin{cor}
 Let $G$ be a graph of order $n$ with a clique partition $F=\{C_1, \ldots, C_k\}$ and let $|C_i|=s_i^F$ for $1\le i\le k$ such that
 $s_1^F\ge s_2^F \ge \cdots \ge s_k^F$. Then
 \begin{equation}\label{equ8}
 \lambda_k(P_G)\ge -\min_{F}s_1^F,
 \end{equation} where minimum is over all clique partitions $F$ of $G$.
 \end{cor}

 In the case of $k>n$, we have $\lambda_i({\cal R}_F)=0$ for $1+n\le i \le k$ by Theorem \ref{pro2}. Since $s_k^F+\lambda_i(P_G)\le \lambda_i({\cal R}_F)\le s_1^F+\lambda_i(P_G)$, we
 get $-s_1^F\le \lambda_i(P_G)\le -s_k^F<0$. We summarize this in the next result.
 \begin{theorem}
 Let $G$ be a graph of order $n$ and a clique partition $F$ with $|F|=k>n$. Then  \\
 $(i)$\, $-s_1^F\le \lambda_i(P_G)\le -s_k^F$~ for $1+n\le i \le k$.\\
 $(ii)$\, $\nu^-(P_G)\ge k-n$.
 \end{theorem}

 The following result follows from (\ref{eq5}).
 \begin{theorem} Let $G$ be a graph of order $n$ with a clique partition $F=\{C_1, \ldots, C_k\}$ and let $|C_i|=s_i^F$ for $1\le i\le k$ such that
 $s_1^F\ge s_2^F \ge \cdots \ge s_k^F$. If $P_G$ is the corresponding clique
 partition graph of $G$, then for $1\le i \le \nu^-(P_G)$,
 \begin{equation}\label{eq31}
 \lambda_{k-i+1}(P_G)\ge -s_i^F.
 \end{equation} Equality in (\ref{eq31}) holds if $G$ is a $s_1^F$ clique-uniform graph with $\nu^-(P_G)=k-n$.
 \end{theorem}

 \vspace{3mm}

 The following concerns the signless Laplacian eigenvalues of a graph.
 \begin{theorem} Let $G$ be a graph of order $n$ and having a clique partition $F$ with $|F|=k$ and assume  $1\le i \le \min\{n,k\}$.\\
 $(i)$\, If $G$ is a $t$ clique-regular graph, then $q_i(G)-\lambda_i(G)\ge t$.\\
 $(ii)$\, If $G$ is a $s$ clique-uniform graph, then $q_i(G)-\lambda_i(P_G)\ge s$.
 \end{theorem}
 \begin{pf}
 From Section \ref{section:matrices}, the signless Laplacian matrix $Q$ of $G$ satisfies $Q\ge {\cal Q}_F$. This fact with Lemma \ref{lem3} gives $q_i(G)\ge
 \lambda_i({\cal Q}_F)$, where $q_i(G)$ and $\lambda_i({\cal Q}_F)$ are respectively, the $i$th largest signless Laplacian eigenvalue of $G$ and the the $i$th largest eigenvalue of matrix ${\cal Q}_F$.
 Using the above analysis combined with Theorem \ref{pro2} and facts $\lambda_i({\cal Q}_F)=t+\lambda_i(G)$ and $\lambda_i({\cal R}_F)=s+\lambda_i(P_G)$ implies the desired results in $(i)$ and $(ii)$.
 \end{pf}

  \subsection{Applications to energy of graphs and matrices}\label{section:energy}
  In this section, using the theory of vertex-clique incidence matrices of a graph, we introduce a new notion of graph energies, as a
  generalization of the incidence energy and the signless Laplacian energy of the
  graph. Finally, we present new
  upper bounds on energies of a graph, its clique partition graph
  and line graph.

 The energy ${\cal E}(G)$ of the graph $G$ defined in (\ref{bm0})
 has the equivalent expressions as follows \cite{AH5}:
 \begin{equation}\label{energy}{\cal E}(G)=2 \sum_{i=1}^{\nu^+} \lambda_i=2 \sum_{i=1}^{\nu^-}
 -\lambda_{n-i+1}=2\max_{1\le k \le n} \sum_{i=1}^k \lambda_i=2\max_{1\le k \le n} \sum_{i=1}^k -\lambda_{n-i+1}
 \end{equation} where $\nu^+$ and $\nu^-$ are
 respectively the positive and the negative inertia of $G$.
 Nikiforov \cite{Ni-1, Ni-2, Ni-3} proposed a significant extension and generalization of the graph energy
 concept.  The energy of an $r \times s$ matrix $B$ is the summation of its singular values, that is, \begin{equation}\label{equ3}{\cal E}(B)=\sum_{i=1}^s \sigma_i(B).\end{equation}

 Consonni and Todeschini \cite{ct} introduced an entire class of
 matrix-based quantities, defined as \begin{equation}\label{equ4}\sum_{i=1}^n |x_i-\overline{x}|,\end{equation} where $x_1,\, x_2,\, \ldots,\, x_n$ are the
 eigenvalues of the respective matrix, and $\overline{x}$ is their
 arithmetic mean.

 According to (\ref{equ3}) and (\ref{equ4}),  two types of energies can then be
 defined for any matrix $B$. The incidence energy $IE(G)$ of a graph $G$ is defined to be the energy of the incidence matrix of
 $G$ of the type (\ref{equ3}), i.e., $$IE(G)={\cal E}(M)=\sum_{i=1}^m \sigma_i(M)=\sum_{i=1}^m  \sqrt{\lambda_i(M^T M)}=\sum_{i=1}^n  \sqrt{\lambda_i(M M^T)}=\sum_{i=1}^n  \sqrt{q_i}.$$
 Similarly, {\it the vertex-clique incidence energy $IE_F(G)$ of $G$ associated with the clique partition $F$} is defined as the energy of the vertex-clique incidence matrix ${\cal
 M}_F$, i.e.,
 \begin{eqnarray}
 IE_F(G)={\cal E}({\cal M}_F)&=&\sum_{i=1}^k \sigma_i({\cal M}_F)=\sum_{i=1}^k  \sqrt{\lambda_i({\cal M}_F^T {\cal M}_F)}\nonumber\\[2mm]
                             &=&\sum_{i=1}^n  \sqrt{\lambda_i({\cal M}_F {\cal M}_F^T)}=\sum_{i=1}^n \sqrt{\lambda_i({\cal Q}_F)}.\nonumber
 \end{eqnarray}
 Observe $$Q-{\cal Q_F}=(D+{\cal A})-({\cal D}_F+{\cal A})=D-{\cal D}_F=diag(d_1-t_1^F, d_2-t_2^F, \ldots, d_n-t_n^F)\ge 0.$$ From the above and using Lemma \ref{lem3} we have
 $q_i=\lambda_i(Q)\ge \lambda_i({\cal Q}_F)$ and, consequently,  we
 have $$IE_F(G)=\sum_{i=1}^n  \sqrt{\lambda_i({\cal Q}_F)}\le \sum_{i=1}^n
 \sqrt{q_i}=IE(G)$$  with equality if and only if $F=E$.

 \vspace{5mm}

 Moreover, $$\sum_{i=1}^n \lambda_i({\cal Q}_F)=\sum_{i=1}^n t_i^F,~~~\sum_{i=1}^n \lambda_i^2({\cal Q}_F)=\sum_{i=1}^n ({(t_i^F)}^2+t_i^F).$$
 Applying the fact that the diagonal entries are majorized by the eigenvalues of ${\cal Q}_F$ and by a similar method given in \cite{AH6} it can be shown that
 $$\sum_{i=1}^n \sqrt{\lambda_i({\cal Q}_F)}\le \sum_{i=1}^n \sqrt{t_i^F}.$$

 \vspace{5mm}

 Considering the energy of the matrix ${\cal Q}_F$ of the type (\ref{equ4}) gives
 \begin{equation}\label{eq6}
 {\cal E}({\cal Q}_F)=\sum_{i=1}^n \left|\lambda_i({\cal
 Q}_F)-\overline{t}\right|,
 \end{equation} where $\overline{t}=\frac{\sum_{i=1}^n t_i^F}{n}$.
 The energy ${\cal E}({\cal Q}_F)$ can be viewed as a generalization
 of the signless Laplacian energy $LE^+(G)$ of $G$ which is defined as follows \cite{Ab}:
 $$LE^+(G)={\cal E}(Q)=\sum_{i=1}^n |q_i-\frac{2m}{n}|.$$
 Due to the similarity of the definitions for signless Laplacian energy
 $LE^+(G)$ and ${\cal E}({\cal Q}_F)$ it follows that in most cases, results derived about $LE^+(G)$ can be generalized to ${\cal E}({\cal Q}_F)$.
 For example, from Lemma 2.12  in \cite{AH5}
 for $LE^+(G)$, we obtain the following:
 \begin{equation}\label{equ7}
 {\cal E}({\cal Q}_F)=\max_{1\le j \le n} \left\{2\sum_{i=1}^j \lambda_i({\cal Q}_F)-2j\,\overline{t}\right\}=2\sum_{i=1}^\tau \lambda_i({\cal Q}_F)-2\overline{t}\,\tau,
 \end{equation}
 where $\tau$ is the largest positive integer such that $\lambda_{\tau}({\cal
 Q}_F)>\overline{t}$.

 Using a method similar to the proof of Corollary 5 in \cite{so} for ${\cal Q}_F-\overline{t}\, I={\cal  D}_F-\overline{t}\,I+{\cal A}$,  we have 
 $${\cal E}({\cal Q}_F)-{\cal E}(G)\le \sum_{i=1}^n |t_i^F-\overline{t}|.$$

 In the next result, we show that for a clique-regular graph $G$ associated with a clique partition $F$, ${\cal E}({\cal Q}_F)={\cal E}(G)$.
 \begin{theorem}\label{thm4}
 If $G$ is a clique-regular graph associated with a clique partition $F$, then ${\cal E}({\cal Q}_F)={\cal E}(G)$.
 \end{theorem}
 \begin{pf} Suppose that $G$ is $t$ clique-regular. Then
 \begin{align}
 {\cal E}({\cal Q}_F)=&\sum_{i=1}^n \left|\lambda_i({\cal Q}_F)-\frac{\sum_{i=1}^n t_i^F}{n}\right|=\sum_{i=1}^n
 \left|\lambda_i({\cal Q}_F)-t\right|\nonumber\\
 =&\sum_{i=1}^n \left|\lambda_i({\cal D}_F+{\cal A}(G))-t\right|=\sum_{i=1}^n \left|\lambda_i(tI+{\cal A}(G))-t\right|\nonumber\\
 =&\sum_{i=1}^n \left|\lambda_i(G)\right|={\cal E}(G).\nonumber
 \end{align} \end{pf}

 Note that for any $t$ clique-regular graph $G$, we have
 $IE_F(G)=\sum_{i=1}^n \sqrt{\lambda_i({\cal Q}_F)}=\sum_{i=1}^n \sqrt{t+\lambda_i}$.
 Next we show that for a clique-uniform graph $G$ we have  ${\cal E}({\cal R}_F)={\cal
 E}(P_G)$.
 \begin{theorem}\label{thm5}
 If $G$ is a clique-uniform graph with the clique partition graph $P_G$, then
   ${\cal E}({\cal R}_F)={\cal E}(P_G)$.
 \end{theorem}
 \begin{pf} Suppose that $G$ is a $s$ clique-uniform graph. Then
 \begin{align}
 {\cal E}({\cal R}_F)=&\sum_{i=1}^k \Big|\lambda_i({\cal R}_F)-\frac{\sum_{i=1}^k
 s_i^F}{k}\Big|=\sum_{i=1}^k \left|\lambda_i({\cal R}_F)-s\right|\nonumber\\
  =&\sum_{i=1}^k \left|\lambda_i({\cal S}_F+{\cal A}(P_G))-s\right|=\sum_{i=1}^k \left|\lambda_i(sI+{\cal A}(P_G))-s\right|\nonumber\\
  =&\sum_{i=1}^k \left|\lambda_i(P_G)\right|={\cal E}(P_G).\nonumber
 \end{align}
 \end{pf}
 Note that for any $s$ clique-uniform graph $G$ with the clique partition graph $P_G$, we have
 $$IE_F(G)=\sum_{i=1}^n \sqrt{\lambda_i({\cal Q}_F)}= \sum_{i=1}^k \sqrt{\lambda_i({\cal R}_F)}=\sum_{i=1}^k \sqrt{s+\lambda_i(P_G)}.$$

 In \cite{AH5} Theorem 3.3, a relation between the energy of the line graph ${\cal E}(L_G)$  and the signless Laplacian energy $LE^+(G)$ of $G$ is
 given. In the following, we generalize this result by using the notion of clique partition of a graph and we provide a comparison between the energy of the clique partition graph ${\cal
 E}(P_G)$ of $P_G$ and ${\cal E}({\cal Q}_F)$. For this we need the following lemma, which is obtained from Theorem \ref{pro2} and is a generalization of Lemma \ref{line}.
 \begin{lemma}\label{lem5}
  Let $G$ be an $s$ clique-uniform graph of order $n$ associated with a clique partition $F$ where $|F|=k$. Then  $$\lambda_i({\cal Q}_F)=\lambda_i(P_G)+s,~~~\mbox{for}~i\in \{1, \ldots, \min\{n, k\}\}.$$
 \end{lemma}

 \begin{theorem}
  Let $G$ be an $s$ clique-uniform graph of order $n$ associated with a clique partition $F$ where $|F|=k$.\\
  $(i)$\,  If $k<n$, then ${\cal E}(P_G)\le {\cal E}({\cal  Q}_F)+\frac{2ks}{n}-2s$.\\
  $(ii)$\, If $k>n$, then ${\cal E}(P_G)\ge {\cal E}({\cal  Q}_F)+\frac{2ks}{n}-2s$.\\
  $(iii)$\, If $k=n$, then ${\cal E}(P_G)={\cal E}({\cal  Q}_F)$.
 \end{theorem}
 \begin{pf} $(i)$\,Let $\nu^+=\nu^+(P_G)\le k <n$. By Lemma
 \ref{lem5} we have
 $$\sum_{i=1}^{\nu^+} \lambda_i(P_G)=\sum_{i=1}^{\nu^+} (\lambda_i({\cal Q}_F)-s)=\sum_{i=1}^{\nu^+} \lambda_i({\cal Q}_F)-s \nu^+.$$
 On the other hand, from (\ref{energy}) we have
 \begin{eqnarray}
 {\cal E}(P_G)=2\sum_{i=1}^{\nu+} \lambda_i(P_G) &=& 2\sum_{i=1}^{\nu+} \lambda_i({\cal  Q}_F)-2s \nu^+ -2\nu^+\frac{\sum_{i=1}^n t_i^F}{n}+2\nu^+\frac{\sum_{i=1}^n t_i^F}{n}\nonumber\\[2mm]
                                                  &\le& {\cal E}({\cal Q}_F)-2s \nu^+ +2\nu^+\frac{ks}{n}~~~\mbox{as}~(\ref{equ7})~,\,\sum_{i=1}^n t_i^F=ks\nonumber\\[2mm]
                                                  &=& {\cal E}({\cal Q}_F)+2\nu^+ \left(\frac{ks}{n}-s\right)~~~\mbox{as}\, \nu^+\ge 1,\, k<n\nonumber\\[2mm]
                                                  &\le& {\cal E}({\cal Q}_F)+\frac{2ks}{n}-2s.\nonumber
 \end{eqnarray}

 $(ii)$\,Recall that $\tau$ is the largest positive integer such that $\lambda_{\tau}({\cal Q}_F)\ge
 \overline{t}=\frac{ks}{n}$ and let $\tau< n< k$. Again by Lemma \ref{lem5} we have
 $$\sum_{i=1}^{\tau} \lambda_i({\cal Q}_F)=\sum_{i=1}^{\tau} (\lambda_i(P_G))+s \tau.$$
 On the other hand, by (\ref{equ7}) and Lemma \ref{lem5} we have
 $$ {\cal E}({\cal Q}_F)=2\sum_{i=1}^{\tau} \lambda_i({\cal Q}_F)-\frac{2ks\tau}{n} = 2\sum_{i=1}^{\tau} \lambda_i(P_G)+2s\tau-\frac{2ks\tau}{n}.$$
 From (\ref{energy}) with the above equation we have
 $${\cal E}(P_G)\ge 2\sum_{i=1}^{\tau} \lambda_i(P_G)={\cal E}({\cal Q}_F)+2\tau\left(\frac{ks}{n}-s\right)\ge {\cal E}({\cal Q}_F)+\frac{2ks}{n}-2s.$$

 $(iii)$\, If $k\neq n$, then ${\cal E}(P_G)\neq {\cal E}({\cal
 Q}_F)$ by $(i)$ and $(ii)$, i.e., if ${\cal E}(P_G)={\cal E}({\cal Q}_F)$, then
 $k=n$. It suffices to show that if $k=n$, then ${\cal E}(P_G)= {\cal E}({\cal
 Q}_F)$. Indeed, if $k=n$, then
 $${\cal E}({\cal Q}_F)=\sum_{i=1}^n |\lambda_i({\cal Q}_F)-\frac{\sum_{i=1}^n t_i^F}{n}|=\sum_{i=1}^n |\lambda_i({\cal Q}_F)-\frac{ks}{n}|.$$
 From the fact $k=n$ with Lemma \ref{lem5} we have ${\cal E}({\cal Q}_F)=\sum_{i=1}^n |\lambda_i(P_G)|={\cal E}(P_G).$
 \end{pf}

 In the following, we present a new upper bound for the energy of a graph $G$.
 \begin{theorem}\label{thm7}
 Let $G$ be a graph of order $n$ and the negative inertia $\nu^-=\nu^-(G)$ and let $t_i^F$ be the $i^{th}$  largest clique degree associated with the clique partition $F$,
 for $1\le i\le n$. Then
 \begin{equation}\nonumber
 {\cal E}(G)\le 2\min_{F} \sum_{i=1}^{\nu^-}  t_i^F,
 \end{equation} where the minimum is given over all clique partitions $F$ of $G$. Equality holds if $G$ is a clique-regular graph associated with a minimum clique partition of size $cp(G)=n-\nu^-$.
 \end{theorem}
 \begin{pf}
 From (\ref{energy}) and (\ref{eq3}) we have $${\cal E}(G)=2\sum_{i=1}^{\nu^-}-\lambda_{n-i+1}\le 2\sum_{i=1}^{\nu^-}
 t_i^F,$$ where $t_i^F$ is $i^{th}$ largest clique-degree of $G$ associated with a
 clique partition $F$. Since this upper bound is valid for any clique
 partition of $G$,  we select the optimal value, namely, $\min\limits_{F}2 \sum_{i=1}^{\nu^-}  t_i^F$. The second part of the proof follows directly from Theorem \ref{thm6}.
 \end{pf}

 \begin{theorem}\label{thm8}
 Let $G$ be a graph of order $n$ with the vertex degrees $d_1\ge d_2 \ge \cdots \ge d_n$. Then
 \begin{equation}\nonumber
 {\cal E}(G)\le 2\sum_{i=1}^h  d_i,
 \end{equation} where $h=\min\{\nu^+,\, \nu^-\}$.
 \end{theorem}
 \begin{pf}
 Considering the fact $t_i^F\le d_i$ for $i\in [n]$ along with Theorem \ref{thm7} gives \begin{equation}\label{equ9}{\cal E}(G)\le 2\sum_{i=1}^{\nu^-}
 d_i.\end{equation} On the other hand, the Laplacian matrix $L=D-{\cal A}$ of $G$ is a 
 positive semi-definite matrix, so ${\cal A}\le
 D$. From this with Lemma \ref{lem3} we obtain $\lambda_i\le d_i$
 for $1\le i \le n$. Then ${\cal E}(G)= 2\sum_{i=1}^{\nu^+}
 \lambda_i\le  2\sum_{i=1}^{\nu^+} d_i$. Using the previous inequality with (\ref{equ9})
 completes the proof.
 \end{pf}

 From Theorem \ref{thm7} with (\ref{alpha}) we obtain the following
 upper bound for the energy of $G$:
 $${\cal E}(G)\le 2\sum_{i=1}^{n-\alpha}  t_i^F\le 2\sum_{i=1}^{n-\alpha}  d_i,$$  where $\alpha$
 is the independence number of the graph $G$.
By (\ref{energy}) and (\ref{eq31}) and applying a similar method done for the proof of Theorem \ref{thm7}, we obtain the next result.
 \begin{theorem}\label{thm9}
 Let $G$ be a graph of order $n$ with a clique partition $F=\{C_1, \ldots, C_k\}$ and let $|C_i|=s_i^F$ for $1\le i\le k$ such that
 $s_1^F\ge s_2^F \ge \cdots \ge s_k^F$. For the clique partition graph
 $P_G$ of $G$, we have
 \begin{equation}\nonumber
 {\cal E}(P_G)\le 2 \min_{F}\sum_{i=1}^{\nu^-(P_G)}  s_i^F.
 \end{equation} Equality holds if $G$ is a clique-uniform graph associated with a minimum clique partition of size $cp(G)=n+\nu^-(P_G)$.
 \end{theorem}

 Next, we present an upper bound for the energy ${\cal E}(L_G)$ of the line graph $L_G$ with a full characterization of the  corresponding
 extreme graphs.
 \begin{theorem}
 Let $G$ be a graph with the line graph $L_G$. Then
 \begin{equation}\label{equ10}
 {\cal E}(L_G)\le 4\,\nu^-(L_G).
 \end{equation} Equality holds if and only if $G$ is a graph with connected components $G_i= (V_i, E_i)$ for $i\ge 1$ with $n_i=|V_i|$ and  $|E_i|\ge 2$, and 
  possibly some isolated vertices or single edges. Further, each non-bipartite connected component $G_i$ satisfies $|E_i|>|V_i|$ and $q_{n_i}\ge
 2$, and each bipartite connected component $G_i$ is either a 4-cycle or satisfies $|E_i|>|V_i|$ and $q_{n_i-1}\ge 2$.
 \end{theorem}
 \begin{pf}
 As previously noted, if the clique partition $F$ of $G$ is as same as the edge set $E$ of $G$,
 then $s_i^F=2$ for $i\in [n]$ and $P_G\cong L_G$. This using Theorem
 \ref{thm9}, we have \begin{equation}\label{equ11}
 {\cal E}(L_G)=2\sum_{i=1}^{\nu^-(L_G)}-\lambda_{m-i+1}(P_G)\le 2 \sum_{i=1}^{\nu^-(L_G)} 2=4\nu^-(L_G),
 \end{equation} which gives the desired result in (\ref{equ10}).

 To characterize these extreme graphs in (\ref{equ10}), we assume equality holds in (\ref{equ11}). Then all negative eigenvalues of $P_G$ must be $-2$ by (\ref{equ11}). We then consider the following two cases:

 $Case ~1)$\, $G$ is connected. First, assume that $m>n$. If $G$ is non-bipartite,
 then by Lemma \ref{line}, $\lambda_i(L_G)=-2$ for $n+1\le i \le m$ and $\lambda_n(L_G)=q_n-2\neq
 -2$ as $q_n\neq 0$. Since $\lambda_n(L_G)$ must be nonnegative, we have $q_n\ge 2$.  Otherwise $G$ is bipartite and by Lemma \ref{line} along with the fact $q_n(G)=0$, $\lambda_i(L_G)=-2$ for $n\le i \le m$ and $\lambda_{n-1}(L_G)=q_{n-1}-2\neq
 -2$ as $q_{n-1}\neq 0$. Since $\lambda_{n-1}(L_G)$ must be nonnegative, it follows that $q_{n-1}\ge 2$. Next, assume that $m=n$. Since all negative eigenvalues of $L_G$ are equal to $-2$, we have $\lambda_m(L_G)=\lambda_n(L_G)=-2$. If $\nu^-=1$,
 then $\Spec(L_G)=\{2, 0, 0, -2\}$ and then $L_G$ is the cycle graph
 $C_4$ of order $4$. Otherwise $\nu^-\ge 2$, and $\lambda_{n-1}=-2$,
 that is, $q_{n-1}=0$, which is a contradiction as $G$ is connected.
 Finally, assume that $m<n$. Since $G$ is connected it must be a tree and
 hence $m=n-1$. In this case we have $\lambda_m(L_G)=\lambda_{n-1}(L_G)=-2$, that
 is, $q_{n-1}=0$, which again leads to a contradiction.

 $Case ~2)$\,Assume  $G$ is disconnected. Since isolated vertices and
 single edges do not affect the negative inertia of $L_G$, we may
 assume that $G$ has connected components along with the possibility of some isolated vertices and
 single edges. Now each connected component of $G$ can be characterized by the first case, and the proof is complete.
 \end{pf}

 \section{Vertex-clique incidence matrix of a graph associated with an edge clique cover}\label{Vertex-clique incidence matrix}

In this section, we consider a slightly more general object of the vertex-clique incidence matrix, denoted by $M_F$,  associated with an edge clique cover $F$ of a graph $G$. Recall that the $(i,j)$-entry of $M_F$ is real and nonzero if and only if the vertex $i$ belongs to the clique $C_j\in F$. 
 To ensure $M_F\,M^T_F\in S(G)$, we arrange the entries of $M_F$ such that the inner product of row $i$ and column $j$ (when $i \neq j$) in $M_F$ is nonzero if and only if $ij\in E(G)$.

As noted in the introduction studying the graph parameter $q(G)$ represents a critical step in the much more general investigation of the inverse eigenvalue problem for graphs. One strategy for minimizing the number of distinct eigenvalues of $M_F\,M_F^T$, we instead consider minimizing the number of distinct eigenvalues of $M_F^T\,M_F$. Consequently achieving an upper bound on the parameter $q(G)$. The key technique used here is to generalize the vertex-clique incidence matrix obtained from an edge clique cover by considering arbitrary positive real entries for $M_F$ or any negative real entries for $M_F$ but paying careful attention to preserving the condition that $M_F\,M_F^T\in S(G)$.

 \subsection{Applications to the minimum distinct eigenvalues of a graph}\label{q=2}
In this section, applying the tool of the vertex-clique incidence matrix of a graph associated with its edge clique cover, we characterize a few new classes of graphs with $q(G)=2$.

If $G$ and $H$ are graphs then the Cartesian product of $G$ and $H$ denoted by $G\square H$, is the graph on the vertex set $V (G)\times V (H)$ with $\{g_1, h_1\}$ and $\{g_2, h_2\}$ adjacent if and only if either $g_1 = g_2$ and $h_1$ and $h_2$ are adjacent in $H$ or $g_1$ and $g_2$ are adjacent in $G$ and $h_1 = h_2$. The first statement in the next theorem can also be found in \cite{AACF}, however, we include a proof here to aid in establishing the second claim in the result below.

\begin{theorem}\label{carGH}
Let $G\cong K_s\square K_2$ with $s\ge 3$. Then $q(G)=2$ and $G$ has an SSP matrix realization with two distinct eigenvalues. 
\end{theorem}
\begin{pf} Let 
 $M=\left(\begin{array}{@{}c@{}}
    M_1 \\
    M_2
  \end{array}\right),$ where $M_1=J_s-(s-1)I_s$ and $M_2=J_s-I_s$. Then we have 
  \begin{equation}\label{blocks}
  A=MM^T=\left(\begin{array}{@{}c@{}}
    M_1 \\
    M_2
  \end{array}\right)\left(M_1^T \, M_2^T\right)=\left(\begin{array}{@{}c|c@{}}
    M_1M_1^T & M_1M_2^T \\\hline
    M_2M_1^T & M_2M_2^T
  \end{array}\right)=  
  \left(\begin{array}{@{}c|c@{}}
    A_1 & (s-1)I_s \\\hline
    (s-1)I_s & A_2
  \end{array}\right),
  \end{equation} where 
  $$A_1=M_1M_1^T=M_1^2=(s-1)^2 I_s+(2-s)J_s\,,~~~ A_2=M_2M_2^T=M_2^2=I_s+(s-2)J_s\,,$$
  $$M_1M_2^T=M_1M_2=(s-1)I_s.$$ From the structure of $A$, we have
  $A\in S(G)$. On the other hand, $$M^TM=\left(M_1^T \, M_2^T\right) \left(\begin{array}{@{}c@{}}
    M_1 \\
    M_2
  \end{array}\right)=M_1^TM_1+M_2^T M_2=(2-s)J_s+(s-1)^2I_s+I_s+(s-2)J_s=cI_s,$$ where $c=s^2-2s+2$. Hence $\Spec(MM^T)=\{c^{[s]},\, 0^{[s]}\}$ and  $q(G)=2$.  
  
  Now, we show that the matrix $A$ has SSP. We need to prove that the only symmetric matrix satisfying  $A\circ X=O$, $I\circ X=O$, and $[A,\, X]=AX-XA=O$ is $X=O$.

From the two equations $A\circ X=O$, $I\circ X=O$, $X$ must have the following form:
$X=\left(\begin{array}{@{}c|c@{}}
    O & X_1 \\\hline
    X_1^T & O
  \end{array}\right)$, where $X_1=\left(
  \begin{array}{cccc}
    0 & x_{12} & \ldots & x_{1s} \\
    x_{21} & 0 &  & x_{2s} \\
    \vdots & \ddots & \ddots & \vdots \\
    x_{s1} & x_{s2} & \ldots & 0 \\
  \end{array}
\right)$. The equality $AX=XA$ gives $X_1=X_1^T$. Also, we have $A_1X_1=X_1A_2$, i.e., $[(s-1)^2 I_s+(2-s)J_s]X_1=X_1 [I_s+(s-2)J_s]$. Hence $sX_1=X_1J_s+J_sX_1$. Then $(sX_1)_{ij}=(X_1J_s+J_sX_1)_{ij}$ for $i,j\in [s]$. Considering $i=j=1$, we have $(sX_1)_{11}=0$ and $(X_1J_s+J_sX_1)_{ii}=2\sum_{j=1}^s x_{1j}$, and then $\sum_{j=1}^s x_{1j}=0$. Considering $(i,j)=(k,k)$ for $2\le k \le s$ we arrive at $\sum_{j=1}^s x_{kj}=0$ for $2\le k \le s$. This means that the row and column sums in $X_1$ are equal to zero. Now,  consider $i, j\in [s]$ where $i \neq j$. We have
$$sx_{ij}=(sX_1)_{ij}=(X_1J_s+J_sX_1)_{ij}=(X_1J_s)_{ij}+(J_sX_1)_{ij}=\sum_{k=1}^s x_{ik}+\sum_{k=1}^s x_{jk}=0.$$ Thus $X_1=O_s$ and consequently, $X=O$. Hence the proof is complete.
\end{pf}

\begin{cor}
 For even $n$, we have $q(\overline{C_n})=2$.
\end{cor}
\begin{pf}
Let $G\cong K_n\backslash H$ and let $H$ be the graph obtained from the complete bipartite graph $K_{n/2, n/2}$ by removing a perfect matching. Then by Theorem \ref{carGH} and Lemma \ref{ssp}, for $H$ or any subgraph of $H$, $q(G)=2$. Considering this with the fact that $C_n$ is a subgraph of $H$, the result is obtained. 
\end{pf}

\begin{theorem}
Let $G$ be a graph obtained from $(K_s\square K_2)\vee sK_1$ by removing a perfect matching between $sK_1$ and a copy of $K_s$. Then $q(G)=2$ and $G$ has an SSP matrix realization with two distinct eigenvalues. 
\end{theorem}
\begin{pf}
Let 
 $M=\left(\begin{array}{@{}c@{}}
    M_1 \\
    M_2\\
    I_s
  \end{array}\right),$ where $M_1=J_s-(s-1)I_s$ and $M_2=J_s-I_s$. Considering the fact that $M_1$ and $M_2$ are symmetric, we have 
  $$A=MM^T=\left(\begin{array}{@{}c@{}}
    M_1 \\
    M_2\\
    I_s
  \end{array}\right)\left(M_1^T \, M_2^T\, I_s \right)=\left(\begin{array}{@{}c|c|c@{}}
    M_1M_1^T & M_1M_2^T & M_1I_s \\\hline
    M_2M_1^T & M_2M_2^T & M_2I_s \\\hline
    M_1  &  M_2  &  I_s
  \end{array}\right)=  
  \left(\begin{array}{@{}c|c|c@{}}
    A_1 & (s-1)I_s &  M_1 \\\hline
   (s-1)I_s & A_2  &  M_2\\\hline
    M_1  &  M_2  &  I_s
  \end{array}\right),$$ where 
  $$A_1=M_1M_1^T=M_1^2=(s-1)^2 I_s+(2-s)J_s\,,~~~ A_2=M_2M_2^T=M_2^2=I_s+(s-2)J_s\,,$$
  $$M_1M_2^T=M_1M_2=(s-1)I_s.$$ From the structure of $A$, we have
  $A\in S(G)$. On the other hand, $$M^TM=\left(M_1^T \, M_2^T \, I_s\right) \left(\begin{array}{@{}c@{}}
    M_1 \\
    M_2 \\
    I_s
  \end{array}\right)=M_1^TM_1+M_2^T M_2+I_s^2=(2-s)J_s+(s-1)^2I_s+I_s+(s-2)J_s+I_s=cI_s,$$ where $c=s^2-2s+3$. This gives $\Spec(MM^T)=\{c^{[s]},\, 0^{[2s]}\}$, which proves $q(G)=2$.  
  
  Now, we show that the matrix $A$ has SSP. We need to prove that the only symmetric matrix satisfying  $A\circ X=O$, $I\circ X=O$, and $[A,\, X]=AX-XA=O$ is $X=O$.

From the two equations $A\circ X=O$, $I\circ X=O$, $X$ must have the following form:
$X=\left(\begin{array}{@{}c|c|c@{}}
        O & X_1 & O \\\hline
    X_1^T & O   & X_2  \\\hline
    O & X_2 &  X_3
  \end{array}\right)$, where $X_1=\left(
  \begin{array}{cccc}
    0 & x_{12} & \ldots & x_{1s} \\
    x_{21} & 0 &  & x_{2s} \\
    \vdots & \ddots & \ddots & \vdots \\
    x_{s1} & x_{s2} & \ldots & 0 \\
  \end{array}
\right)$, $X_2=diag(y_1, \ldots, y_s)$\, and \,$X_3=\left(
  \begin{array}{cccc}
    0 & z_{12} & \ldots & z_{1s} \\
    z_{12} & 0 &  & z_{2s} \\
    \vdots & \ddots & \ddots & \vdots \\
    z_{1s} & z_{2s} & \ldots & 0 \\
  \end{array}
\right)$. The matrix equation
\begin{equation}
AX=XA  \label{matequ1}
\end{equation} gives $X_1=X_1^T$. From (\ref{matequ1}) we also have
$M_2X_2+X_3=X_2M_2+X_3$, i.e., $(J_s-I_s)X_2=X_2 (J_s-I_s)$, i.e., $J_sX_2=X_2J_s$. This gives  $y_1=y_2=\cdots=y_s$, i.e., $X_2=y_1I_s$.

Again from (\ref{matequ1}), we have $A_1X_1+M_1X_2=X_1A_2$, that is, $M_1X_2=X_1A_2-A_1X_1$, that is, 
$(J_s-(s-1)I_s)(y_1I_s)=X_1 (I_s+(s-2)J_s)-((s-1)^2 I_s+(2-s)J_s)X_1$, i.e., 
$$y_1(2-s)I_s+y_1J_s=(2s-s^2)X_1+(s-2)X_1J_s+(s-2)J_sX_1.$$ Considering a main diagonal entry, say $(i,i)$,  in the above matrix equation, we obtain 
\begin{equation}
\sum_{j=1}^s x_{ij}=-\frac{y_1}{2}.  \label{2}
\end{equation} 
Considering the $(i,j)$-entry  in the above matrix equation, we obtain
$x_{ij}=-y_1\frac{s-1}{s-2}$. From the above and (\ref{2}), $y_1=0$, that is, $X_2=O$. Using the equation $A_1X_1+M_1X_2=X_1A_2$, we arrive at the matrix equation $A_1X_1=X_1A_2$. Following a similar argument as in the proof of Theorem \ref{carGH} we obtain $X_1=O$. 

Again from (\ref{matequ1}), we have $M_1X_1+X_2=X_2A_2+X_3M_2$. Since $X_1=X_2=O$, we get 
$X_3M_2=O$, i.e. $X_3=X_3J_s$. Considering both the $(i,i)$ and $(i,j)$ entries from the matrix equation, we arrive at $\sum_{k=1}^s z_{ik}=0$ and $z_{ij}=\sum_{k=1}^s z_{ik}=0$, that is, $X_3=O$, which  gives $X=O$. 
\end{pf}

\begin{cor}
 Consider the complete bipartite graph $K_{s,s}$ by removing a perfect matching. Define a new graph $H$ by adding a copy of $K_s$ to this graph such that each vertex in $K_s$ is adjacent to the corresponding vertex in a copy of $sK_1$. Then $q(\overline{H})=2$. Moreover, the result holds for any subgraph of $H$ on the same vertex set. 
\end{cor}

In \cite{RT}, the authors studied the problem of graphs requiring property $p(r,s)$. A graph $G$ has $p(r,s)$ if it contains a path of length $r$ and every path of length $r$ is contained in a cycle of length $s$. They prove that the smallest integer $m$ so that every graph on $n$ vertices with   $m$ edges has $p(2,4)$ (or each path of length $2$ is contained either  in a $3$-cycle, or a $4$-cycle) is ${n\choose 2}-(n-4)$ for all $n\ge 5$.
Using this, it was noted in \cite{private} that the above equation from \cite{RT} implies that the fewest number of edges required to guarantee that all graphs $G$ on $n$ vertices satisfy $q(G)=2$ is at least  ${n\choose 2}-(n-3)$. For small values of $n$, it is known that in fact, equality holds in the previous claim. Namely, if at most $n-3$ edges are removed from the complete graph $K_n$ with $n \leq 7$, then the resulting graph has a matrix realization with two distinct eigenvalues.
Along these lines and based on \cite{private}  the following is a natural conjecture:

\begin{conj}\label{conj1}
Removing up to $n-3$ edges from $K_n$ does not change the number of distinct eigenvalues of $K_n$. That is,  for any subgraph $H$ of $K_n$ with $|E(H)|\le n-3$
$$q(K_n\backslash H)=2.$$
\end{conj}

We confirm Conjecture \ref{conj1} for $n=7,8$ and note that our analysis of the case $n=7$ differs slightly from \cite{private}. For this, we need the next few lemmas.

\begin{lemma}\label{lemT1}
 Let $T_1$ be the tree given in Figure \ref{tree}. We have $q(\overline{T_1})=2$ and $\overline{T_1}$ has an SSP matrix realization with two distinct eigenvalues.
 \begin{figure}[htb]
\begin{center}
 \includegraphics[height=4cm,keepaspectratio]{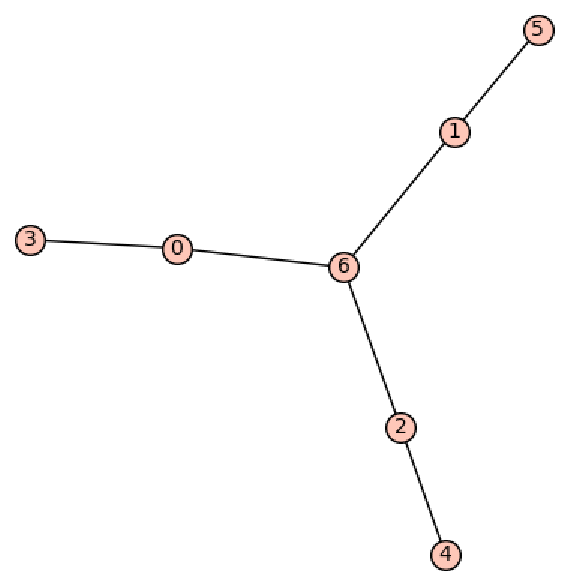}\\
 \caption{Tree $T_1$.}\label{tree} 
 \end{center}
 \end{figure}
\end{lemma}
\begin{pf}
Consider the $7\times 4$ matrix $M_1$ as follows:
$$M_1=\left(
  \begin{array}{cccc}
     1 & -2 & 2 & 1  \\
     2 & -1 & -2 & 2 \\
     2 & 2  & 1 & 2 \\
     1 & 2  & 2 & 0  \\
    -2 & -1 & 2 & 0  \\
     2 & -2 & 1 & 0 \\
     1 & 0  & 0 & 0 \\
  \end{array}
\right).$$ Using the Gram-Schmidt method we can arrive at a column orthonormal matrix $M_2$. In this case we have
$A=M_2M_2^T\in S(\overline{T_1})$. Also $M_2^TM_2=I_4$ and then $\Spec(A)=\{1^{[4]},\, 0^{[3]}\}$. This proves that $q(\overline{T_1})=2$. Furthermore, $A$ has SSP (this can be confirmed using SageMath), and  by Lemma \ref{ssp}, the complement of any subgraph of $T_1$ on the same vertex set also has a matrix realization with two distinct eigenvalues.  
\end{pf}

\begin{figure}[htb]
\begin{center}
 \includegraphics[height=3.8cm,keepaspectratio]{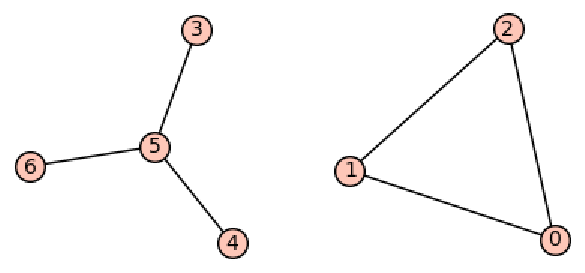} 
 \caption{The graph $G$.} 
 \label{fig 4.7}
 \end{center}
 \end{figure}

\begin{lemma}\label{lemSC}
 Let $G\cong K_{1,3} \cup K_3$. Then $q(\overline{G})=2$ and $\overline{G}$ has an SSP matrix realization with two distinct eigenvalues.
\end{lemma}
\begin{pf}
Consider the $7\times 3$ matrix $M_1$ corresponding to the labeled graph $G$ given in Figure \ref{fig 4.7} as follows:
$$M_1=\left(
  \begin{array}{cccc}
     1 & 2   & 2   \\
     2 & 1   & -2  \\
     2 & -2  & 1  \\
     1 & 1   & 1   \\
     1 & -1  & 1  \\
-\sqrt{2}& 0 & \sqrt{2}  \\
     0 & \sqrt{2}  & 0  \\
  \end{array}
\right).$$ 
$A=M_1M_1^T\in S(\overline{G})$. Also $M_1^TM_1=13I_3$ and then $\Spec(A)=\{13^{[3]},\, 0^{[4]}\}$. This proves that $q(\overline{G})=2$. Furthermore, $A$ has SSP (a computation that can be verified by SageMath), and by Lemma \ref{ssp}, the complement of any subgraph of $G$ on the same vertex set also has a matrix realization with two distinct eigenvalues.  
\end{pf}

We now verify that Conjecture \ref{conj1} holds for $n=7$.
\begin{theorem}
Removing up to $4$ edges from $K_7$ does not change the number of distinct eigenvalues of $K_7$, i.e., for any subgraph $H$ of $K_7$ on 7 vertices, with $|E(H)|\le 4$ we have 
$$q(K_7\backslash H)=2.$$
\end{theorem}
\begin{pf} It suffices to show $\overline{H}$ for any graph $H$ in Figure \ref{G7-4} has a matrix realization with two distinct eigenvalues. Suppose that the graphs in Figure \ref{G7-4} are denoted by $H_i$  for $i\in [10]$ from left to right in each row. Then the graphs $H_i$ for $i=1,3,7,8,10$ are the union of complete bipartite graphs with some isolated vertices. By Lemma \ref{MB} (2), the complements of these graphs and any subgraph of these graphs have a matrix realization with two distinct eigenvalues. 
Also $q(\overline{H_i})=2$ for $i=4, 5, 9$ and for any subgraph $H'_i$ of $H_i$, $q(\overline{H'_i})=2$ by Lemma \ref{lemT1}.
Moreover, $q(\overline{H_6})=2$ and for any subgraph $H'_6$ of $H_6$, $q(\overline{H'_6})=2$  by Lemma \ref{lemSC}. Additionally, from Lemmas \ref{lemT1} and \ref{lemSC} such realizations exists with the SSP. Hence any subgraph of these graphs have a matrix realization with two distinct eigenvalues. 
To complete the proof, we only need to show the complement graph of $H_2$ has a matrix realization with two distinct eigenvalues with the SSP. To this end, consider the $7\times 3$ matrix $M_1$ as follows:
$$M_1=\left(
  \begin{array}{ccc}
     1 & -2 & 1  \\
     2 & -1 & 2 \\
     2 & 2  & 2 \\
     1 & 2  & 0  \\
    -2 & -1 & 0  \\
     2 & -2 & 0 \\
     1 & 0  & 0 \\
  \end{array}
\right).$$ Using the Gram-Schmidt method we can arrive at a column orthonormal matrix $M_2$. We have
$A=M_2M_2^T\in S(\overline{H_2})$. Also $M_2^TM_2=I_3$ and then $\Spec(A)=\{1^{[3]},\, 0^{[4]}\}$. Hence $q(\overline{H_2})=2$. Furthermore, $A$ has SSP (a computation that van be verified by SageMath), and by Lemma \ref{ssp}, the complement of any subgraph of $H_2$ on the same vertex set also has a matrix realization with two distinct. 
\begin{figure}[htb]
\begin{center}
 \includegraphics[height=8cm,keepaspectratio]{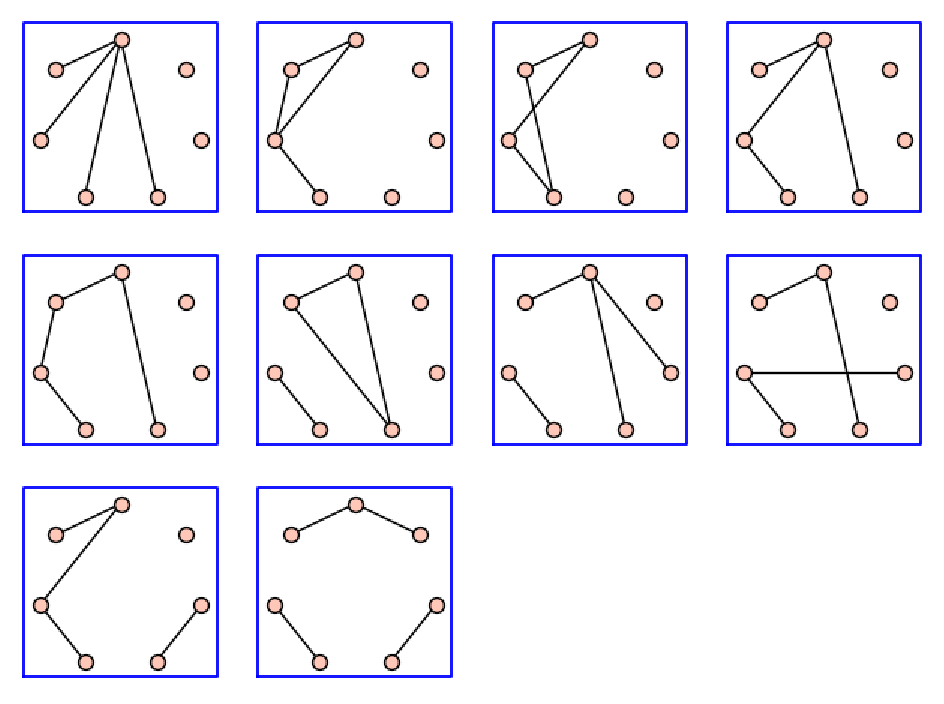}
 \vspace*{-3mm}
 \caption{All graphs with 7 vertices and 4 edges.}
 \label{G7-4}
 \end{center}
 \end{figure}
\end{pf} \vspace*{-1cm}

We require the following results to confirm Conjecture \ref{conj1} for $n=8$. 

\begin{lemma}\label{forG4}
 Let $G\cong H_1\cup 2K_1$, where $H_1$ is the graph on the left given in Figure \ref{L4.9}. Then $q(\overline{G})=2$ and $\overline{G}$ has an SSP matrix realization with two distinct eigenvalues.
\end{lemma}
\begin{pf}
Given $G$ as assumed it can be shown without too much difficulty that  $\overline{G}\cong (H_2 \vee K_3)-e$, where $H_2$ is the graph on the right given in Figure \ref{L4.9} minus an edge $e$ with one endpoint in $K_3$ and the other endpoint in $H_2$ with degree three. 
\begin{figure}[htb]
\begin{center}
\hspace*{4.5cm} \includegraphics[height=7cm,keepaspectratio]{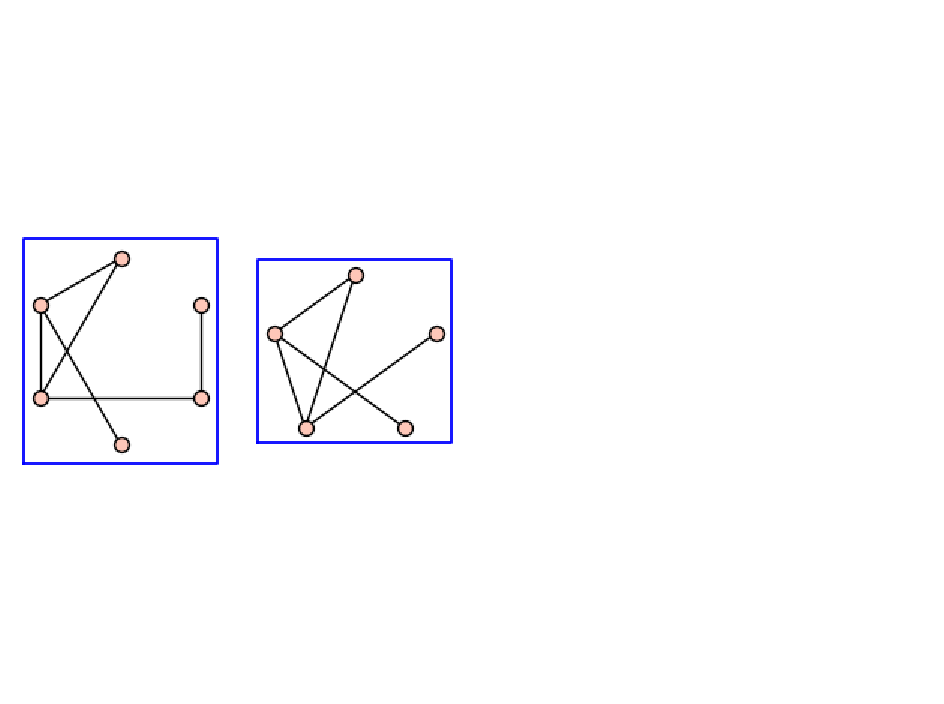}
\vspace*{-2cm}
 \caption{The graphs $H_1$ (left) and $H_2$ (right).}
 \label{L4.9}
 \end{center}
\end{figure}
 Suppose $M=\left(\begin{array}{@{}c@{}}
    M_1 \\
    M_2
  \end{array}\right),$ is a vertex-clique incidence matrix of $\overline{G}$, where the blocks $M_1$ and $M_2$ are vertex-clique incidence matrices corresponding to graphs $H_2$ and $K_3$, that is, $MM^T\in S(\overline{G}).$ From (\ref{blocks}) we have  $M_1M_1^T\in S(H_2)$ and $M_2M_2^T\in S(K_3)$. On the other hand, we have 
  \begin{equation}\label{tblocks}
  M^TM=M_1^TM_1+M_2^T M_2.
  \end{equation} Consider a vertex-clique incidence matrix $M_1$ as follows:
$$M_1=\left(
  \begin{array}{ccc}
     1 & 0 & 0  \\
     1 & 0 & 1 \\
     1 & 1  & 0 \\
     0 & \sqrt{2}  & 0  \\
     0 & 0 & \sqrt{2}  
  \end{array}
\right).$$ Then we have $M_1M_1^T\in S(H_2)$ and $M_1^TM_1=\left(
  \begin{array}{ccc}
     3 & 1 & 1  \\
     1 & 3 & 0 \\
     1 & 0  & 3
  \end{array}
\right)$. Given $M_1$ above, the remainder of the proof is devoted to constructing a matrix $M_2$ so that following  (\ref{tblocks}) we have $M^TM=cI_3$, for some scalar $c$. Consider a matrix $M_2$ so that  
\begin{equation}\label{M20}
M_2^TM_2=\left(
  \begin{array}{ccc}
     a & -1 & -1  \\
     -1 & a & 0 \\
     -1 & 0  & a
  \end{array}
\right),\end{equation} where $a$ is a constant. Suppose the matrix $M_2=\left(
  \begin{array}{ccc}
     x_1 & y_1 & z_1  \\
     x_2 & y_2 & z_2 \\
     x_3 & y_3  & z_3
  \end{array}
\right).$ This with (\ref{M20}) leads to the following equations:
$$x_1^2+x_2^2+x_3^2=y_1^2+y_2^2+y_3^2=z_1^2+z_2^2+z_3^2=a,$$ 
$$x_1y_1+x_2y_2+x_3y_3=-1,~~x_1z_1+x_2z_2+x_3z_3=-1,~~y_1z_1+y_2z_2+y_3z_3=0.$$
Solving this system of non-linear equations, we have a candidate matrix $M_2$:
 $M_2=\left(
  \begin{array}{ccc}
     1 & -1 & z_1  \\
     -1 & 2 & z_2 \\
     2 & 1  & z_3
  \end{array}
\right)$, where $z_1=\frac{1}{7}(2\sqrt{51}-1)$, $z_2=\frac{1}{35}(6\sqrt{51}+4)$, and $z_3=\frac{-1}{35}(2\sqrt{51}+13)$. Thus
$$M=\left(
  \begin{array}{ccc}
     1 & 0 & 0  \\
     1 & 0 & 1 \\
     1 & 1  & 0 \\
     0 & \sqrt{2}  & 0  \\
     0 & 0 & \sqrt{2}  \\
     \hline
     1 & -1 & z_1 \\
    -1 & 2  & z_2 \\
     2 & 1  & z_3 
  \end{array}
\right).$$ It is obvious that $MM^T\in S(\overline{G})$ and $M^TM=9I_3$. Then by the fact that matrices $AB$ and $BA$ have same nonzero eigenvalues, we have $\Spec(MM^T)=\{9^{[3]},\, 0^{[5]}\}$, and then $q(\overline{G})=2$. Moreover, applying a basic computation from SageMath, we can confirm that $MM^T$ has SSP and this completes the proof.
\end{pf}

By Lemma \ref{forG4}, $\overline{G}$ has an SSP realization $A=MM^T$ with two distinct eigenvalues. Then by Lemma \ref{ssp}, any supergraph on the same vertex set as $G$ has a realization with the same spectrum as $A$. In particular, $q(H_2 \vee K_3)=2$. 
This is stated in the following corollary. 

\begin{cor}\label{cor-G4}
 Let $G\cong H_2\cup 3K_1$, where $H_2$ is the right graph given in Figure \ref{L4.9}. Then $q(\overline{G})=2$ and $\overline{G}$ has an SSP matrix realization with two distinct eigenvalues.
\end{cor}

\begin{lemma}\label{forG6}
 Let $G\cong H_3\cup 3K_1$, where $H_3$ is obtained from $C_5$ by joining a vertex to any vertex in $C_5$. Then $q(\overline{G})=2$ and $\overline{G}$ has an SSP matrix realization with two distinct eigenvalues.
\end{lemma}
\begin{pf}
We know that $\overline{G}\cong (C_5\vee K_3)-e$, where $e$ is an edge with one endpoint in $K_3$ and the other in $C_5$. 
 Suppose $M=\left(\begin{array}{@{}c@{}}
    M_1 \\
    M_2
  \end{array}\right),$ is a vertex-clique incidence matrix of $\overline{G}$, where blocks $M_1$ and $M_2$ are vertex-clique incidence matrices corresponding to graphs $C_5$ and $K_3$, that is, $MM^T\in S(\overline{G}).$ From (\ref{blocks}) we have  $M_1M_1^T\in S(C_5)$ and $M_2M_2^T\in S(K_3)$. On the other hand, we also have the equations in  (\ref{tblocks}). Now, we consider a vertex-clique incidence matrix $M_1$ as follows:
$$M_1=\left(
  \begin{array}{ccc}
     1 & 0 & 0  \\
     1 & 1 & 0 \\
     -1 & 1  & 1 \\
     0 & -1  & 1  \\
     0 & 0 & 1  
  \end{array}
\right).$$ Then $M_1M_1^T\in S(C_5)$ and $M_1^TM_1=\left(
  \begin{array}{ccc}
     3 & 0 & -1  \\
     0 & 3 & 0 \\
     -1 & 0  & 3
  \end{array}
\right)$.  Given $M_1$ above, the remainder of the proof is devoted to constructing a matrix $M_2$ so that following  (\ref{tblocks}) we have $M^TM=cI_3$, for some scalar $c$. We need to create a matrix $M_2$ so that 
\begin{equation}\label{M2}
M_2^TM_2=\left(
  \begin{array}{ccc}
     a & 0 & 1  \\
     0 & a & 0 \\
     1 & 0  & a
  \end{array}
\right),\end{equation} where $a$ is a constant. Suppose  $M_2=\left(
  \begin{array}{ccc}
     x_1 & y_1 & z_1  \\
     x_2 & y_2 & z_2 \\
     x_3 & y_3  & z_3
  \end{array}
\right).$ This with (\ref{M2}) leads to the following equations:
$$x_1^2+x_2^2+x_3^2=y_1^2+y_2^2+y_3^2=z_1^2+z_2^2+z_3^2=a,$$ 
$$x_1y_1+x_2y_2+x_3y_3=0,~~x_1z_1+x_2z_2+x_3z_3=1,~~y_1z_1+y_2z_2+y_3z_3=0.$$
Solving these non-linear equations we have
 $M_2=\left(
  \begin{array}{ccc}
     \frac{1}{\sqrt{3}} & 0 &  \frac{1}{\sqrt{3}}  \\
     \frac{1}{\sqrt{3}} &  \frac{1}{\sqrt{2}} &  \frac{1}{\sqrt{3}} \\
     \frac{1}{\sqrt{3}} &  \frac{-1}{\sqrt{2}}  &  \frac{1}{\sqrt{3}}
  \end{array}
\right)$. Thus we have
$$M=\left(
  \begin{array}{ccc}
     1 & 0 & 0  \\
     1 & 1 & 0 \\
     -1 & 1  & 1 \\
     0 & -1  & 1  \\
     0 & 0 & 1  \\
     \hline
     \frac{1}{\sqrt{3}} & 0 &  \frac{1}{\sqrt{3}}  \\
     \frac{1}{\sqrt{3}} &  \frac{1}{\sqrt{2}} &  \frac{1}{\sqrt{3}} \\
     \frac{1}{\sqrt{3}} &  \frac{-1}{\sqrt{2}}  &  \frac{1}{\sqrt{3}}
  \end{array}
\right).$$ It is clear that $MM^T\in S(\overline{G})$ and $M^TM=4I_3$. Then by the fact that matrices $AB$ and $BA$ have same nonzero eigenvalues, we have $\Spec(MM^T)=\{4^{[3]},\, 0^{[5]}\}$, and $q(\overline{G})=2$. Moreover, applying a basic computation from  SageMath, it follows that $MM^T$ has SSP and this completes the proof.
\end{pf}

By Lemma \ref{forG6}, $\overline{G}$ has an SSP realization $A=MM^T$ with two distinct eigenvalues. By Lemma \ref{ssp}, any supergraph on the same set of vertices as $G$ has a matrix realization with same spectrum as $A$. Thus $q(C_5 \vee K_3)=2$. This is stated in the following corollary. 
\begin{cor}\label{cor-G6}
 Let $G\cong C_5\cup 3K_1$. Then $q(\overline{G})=2$ and $\overline{G}$ has an SSP matrix realization with two distinct eigenvalues.
\end{cor}

\begin{figure}[htb]
\begin{center}
 \includegraphics[height=3.8cm,keepaspectratio]{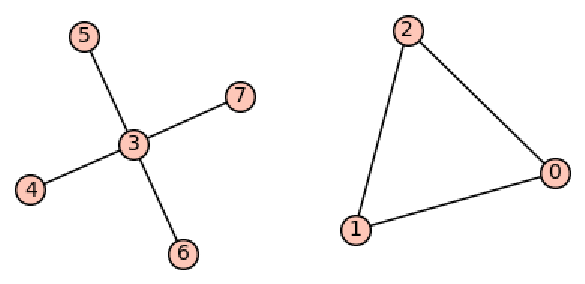} 
 \caption{The graph $G$.} 
 \label{fig 4.13}
 \end{center}
 \end{figure} \vspace{-5mm}

\begin{prop}\label{}
 Let $G\cong K_3\cup K_{1,n-4}$, where $n\ge 7$. Then $q(\overline{G})=2$ and $\overline{G}$ has an SSP matrix realization with two distinct eigenvalues.
\end{prop}
\begin{pf}
We show that the complement of $G$ has a matrix realization with two distinct eigenvalues with the SSP. Consider $n\times 3$ matrix $M_1$ with rows labeled as given in Figure \ref{fig 4.13} for $n=8$:
$$M_1=\left(
  \begin{array}{cccc}
     1 & 2   & 2   \\
     2 & 1   & -2  \\
     2 & -2  & 1  \\
-\sqrt{2}& 0 & \sqrt{2}  \\
     0 & \sqrt{\frac{2}{n-4}}  & 0  \\
     \vdots & \vdots & \vdots \\
     0 & \sqrt{\frac{2}{n-4}}  & 0  
  \end{array}
\right).$$ We have
$A=M_1M_1^T\in S(\overline{G})$. Also $M_1^TM_1=11\,I_3$ and then $\Spec(A)=\{11^{[3]},\, 0^{[n-3]}\}$. This proves that $q(\overline{G})=2$. To verify that $A$ has SSP, suppose $X$ is a symmetric matrix such that  $A\circ X=O$, $I\circ X=O$, and $[A,\, X]=AX-XA=O$. Note to verify $[A,\, X]=AX-XA=O$ it is equivalent to prove that $AX$ is symmetric. Now assume that $X$ has the form:
\[X=\left(\begin{array}{@{}c|c|c@{}}
        0 & O & x^T \\\hline
    O & X_1   & O  \\\hline
    x & O &  O
  \end{array}\right),
  \, {\rm where} \, X_1=\left(
  \begin{array}{ccc}
    0 & a & b \\
    a & 0 & c \\
    b & c & 0 \\
  \end{array}
\right),\] and $x$ is a (possibly) nonzero vector of size $n-4$. Since $AX$ is symmetric, comparing the (1,3) and (3,1) blocks of $AX$ we note that $\alpha Jx = 4x$. So if we set $\beta=\one^T x$, then $x=\frac{\alpha}{4}\beta \one$.
Comparing the (1,2) and (2,1) blocks of $AX$ gives
\[ 2\sqrt{\alpha} \beta = -4\sqrt{2}a - \sqrt{2}b=-\sqrt{2}b + 4\sqrt{2}c, \; {\rm and} \; \sqrt{\alpha}\beta = \sqrt{2}a - \sqrt{2}c.\]
Hence it follows that $a=-c$ and $\beta =\frac{2\sqrt{2}a}{\sqrt{\alpha}}$. Finally, comparing the (2,3) and (3,2) blocks of $AX$, we have 
\[ a\sqrt{\alpha} -2b\sqrt{\alpha} = 2a\sqrt{\alpha} -2c\sqrt{\alpha}=
2b\sqrt{\alpha} +c\sqrt{\alpha} = \left(\frac{\alpha}{4}\beta\right)^2 = \frac{a^2}{2\alpha}. \] From the above equations we deduce that $b=-\frac{3}{2}a$. Substituting the equations $a=-c$, $\beta =\frac{2\sqrt{2}a}{\sqrt{\alpha}}$, and $b=-\frac{3}{2}a$ into the equation $2\sqrt{\alpha} \beta = -\sqrt{2}b + 4\sqrt{2}c$, yields $4\sqrt{2} a = \frac{3}{\sqrt{2}} a -4\sqrt{2}a$. Assuming $a\neq 0$, implies an immediate contradiction. Thus $a=0$, and it follows, based on the analysis above that $X=0$. Hence $A$ has the SSP. Using the fact that this matrix realization has the SSP together with Lemma \ref{ssp}, it follows that the complement of any subgraph of $G$ on the same vertex set also realizes distinct eigenvalues. 
\end{pf}

\begin{figure}[htb]
\begin{center}
 \includegraphics[height=4cm,keepaspectratio]{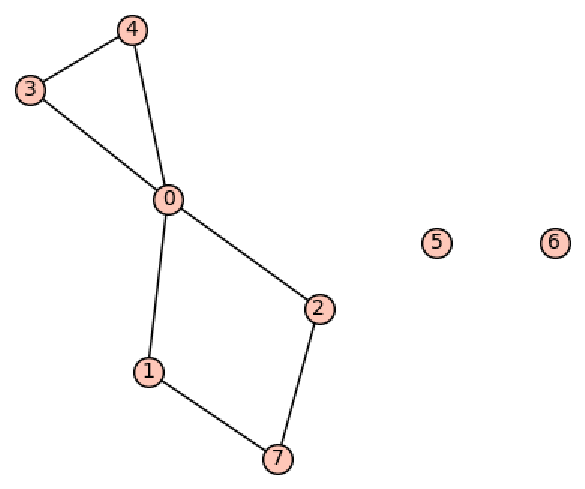} 
 \caption{The graph $G$.} 
 \label{fig 4.14}
 \end{center}
 \end{figure} 
 \vspace*{-1cm}

\begin{lemma}\label{tri-squ}
 Let $G$ be the graph given in Figure \ref{fig 4.14}. Then $q(\overline{G})=2$ and $\overline{G}$ has an SSP matrix realization with two distinct eigenvalues.
\end{lemma}
\begin{pf}
We show that the complement graph of $G$ has a matrix realization with two distinct eigenvalues with the SSP. To do this, first we consider $8\times 3$ matrix $M$ as follows:
$$M=\left(
  \begin{array}{cccc}
     \sqrt{\frac{15}{2}} & 0 & 0   \\
     0 & 1  & 1  \\
     0 & 1  & 1  \\
     0 & 1  & 2   \\
     0 & -2 & 1  \\
     1 & -1 & 0  \\
     1 & 0  & 1  \\
  \sqrt{\frac{2}{2}} & \sqrt{2}  & -\sqrt{2}  
  \end{array}
\right).$$ We have
$A=MM^T\in S(\overline{G})$. Also $M^TM=10\,I_3$ so $\Spec(A)=\{10^{[3]},\, 0^{[5]}\}$. This proves that $q(\overline{G})=2$. Furthermore, $A$ has SSP (observed using SageMath) and by Lemma \ref{ssp}, the complement of any subgraph of $G$ on the same vertex has a matrix realization having two distinct eigenvalues. 
\end{pf}

\begin{figure}[htb]
\begin{center} \includegraphics[height=9cm,keepaspectratio]{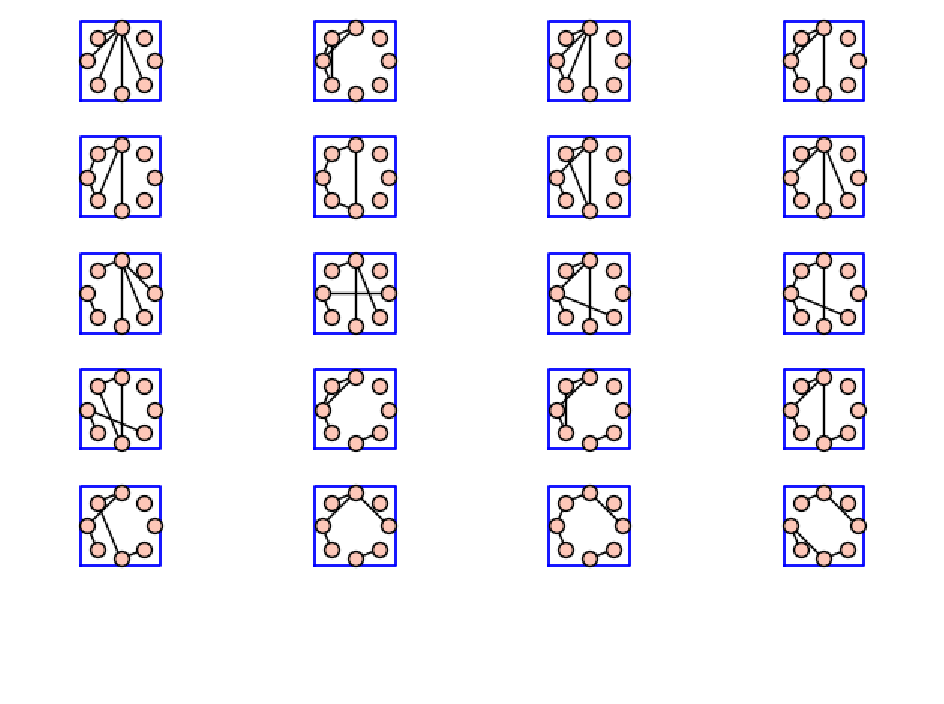} 
 \includegraphics[height=5cm,keepaspectratio]{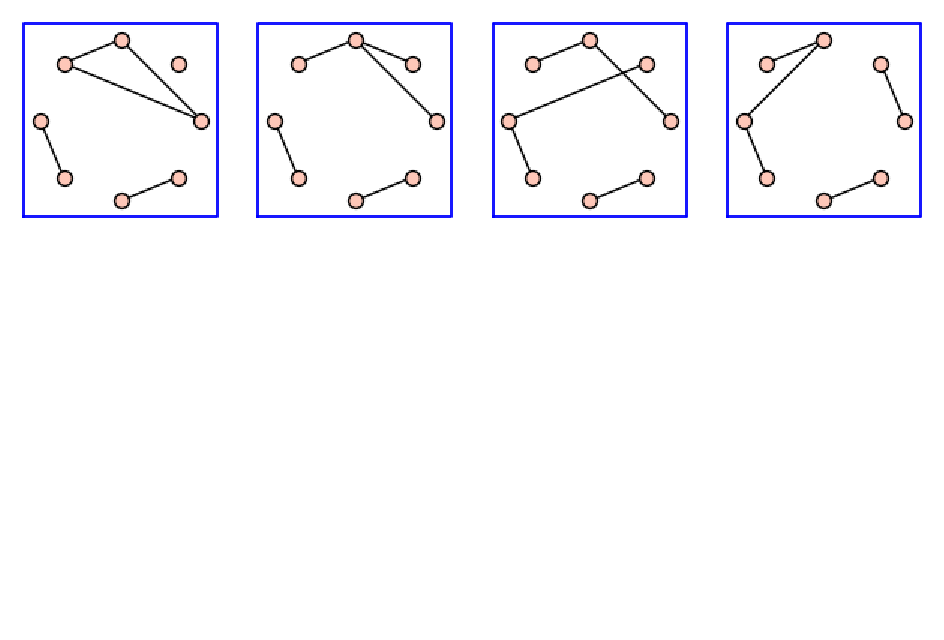} \vspace*{-3cm}
 \caption{All graphs with 8 vertices and 5 edges.} 
 \label{fig}
 \end{center}
 \end{figure}

Now we are in a position to establish that Conjecture \ref{conj1} holds for $n=8$.
\begin{theorem}
Removing up to $5$ edges from $K_8$ does not change the number of distinct eigenvalues of $K_8$, i.e., for any subgraph $H$ on 8 vertices of $K_8$ with $|E(H)|\le 5$,
$$q(K_8\backslash H)=2.$$
\end{theorem}
\begin{pf}
It suffices to show $\overline{H}$ for any graph $H$ in Figure \ref{fig} has a matrix realization with two distinct eigenvalues. Suppose that the graphs in Figure \ref{fig} are denoted by $H_i$  for $i\in [24]$ from left to right in each row. 
The graphs $H_i$ for $i=1,2,9,10,15,22,23$ are the union of complete bipartite graphs with some isolated vertices. By Lemma \ref{MB} (2), the complements of these graphs and any subgraph of these graphs have a matrix realization with two distinct eigenvalues. Also $q(\overline{H_i})=2$ for $i=5, 11,12,16,17,18,19,20,24$ and for any subgraph $H'_i$ of $H_i$, $q(\overline{H'_i})=2$ by Theorem \ref{carGH}. 
For $i=3,7, 8, 13, 14$, we have $q(\overline{H_i})=2$ and for any subgraph $H'_i$ of $H_i$, $q(\overline{H'_i})=2$ by Lemma \ref{tri-squ}. Additionally, from Theorem \ref{carGH} and Lemma \ref{tri-squ} such realizations exists with the SSP. Hence any subgraph of these graphs have a matrix realization with two distinct eigenvalues. 

Further $q(\overline{H_{21}})=q(\overline{(2K_2\cup K_1)\cup K_3})=q(G\vee 3K_1)=2$ by Lemma \ref{lemma-join}, where the graph 
$G=\overline{2K_2\cup K_1}=K_{2,2}\vee K_1$ is connected. If we remove any edges in $H_{21}$ from the triangle, then the complement of the result graph has at least two distinct eigenvalues by Lemma \ref{MB} (2), and if we remove any edges in $H_{21}$ from out of the triangle, again by \ref{lemma-join} we can see that the complement of the result graph has at least two distinct eigenvalues. 
We have $q(\overline{H_4})=2$ and the complement of any subgraph of this graph has a matrix realization with two distinct eigenvalues, by Corollary \ref{cor-G4}.  Moreover,  $q(\overline{H_6})=2$, and the complement of any subgraph of this graph also has a matrix realization with two distinct eigenvalues, by Corollary \ref{cor-G6}. This completes the proof of the theorem.
\end{pf}

 \section{Concluding remarks and open problems}
 In this work, we utilized the notions of a clique partition and an edge clique cover of a graph to introduce and explore the various properties of a vertex-clique incidence matrix of the graph, which can be viewed as a generalization of the vertex-edge incidence matrix.
 Using these new incidence matrices, we obtained sharp interesting lower bounds concerning the 
 negative eigenvalues and thus the negative inertia of a graph, and
 we generalize the notion of the line graph of a graph by introducing the clique partition graph of the given graph. Additionally, we determined 
 the relations between the spectrum of a graph and its clique partition graph.
 Further, we generalized the notion of incidence energy and signless Laplacian energy of a graph and
 provided some novel upper bounds for the energies of a graph, its clique partition graph, and the line graph. Finally, applying a general version of a vertex-clique incidence matrix of a graph associated with its edge clique cover, we were able to characterize a few classes of graphs with $q(G)=2$. To close we list two important and unresolved issues related to some of the content of the current work.

 \noindent \\
 {\bf Problem 1:} Characterize the corresponding extreme graphs for which the inequalities given in (\ref{eq2}), (\ref{equ1}), (\ref{eq3}), (\ref{equ8}), and (\ref{eq31}) hold with equality.

 \noindent \\ 
 {\bf Problem 2:} Prove that Conjecture \ref{conj1} is valid for any graph $G$ of order at least 9. \\

\noindent {\bf Acknowledgements} \\
Dr. Fallat's research was supported in part by an NSERC Discovery Research Grant, Application No.: RGPIN-2019-03934.


\begin{thebibliography}{99}

 \bibitem{Ab} N. Abreu, D. M. Cardoso, I. Gutman, E. A. Martins, M. Robbiano, Bounds for the signless Laplacian energy, Linear Algebra Appl., 435 (2011) 2365--2374.


 \bibitem{AACF} B. Ahmadi, F. Alinaghipour, M.S. Cavers, S. M. Fallat, K. Meagher, and S. Nasserasr, Minimum number of distinct eigenvalues of graphs, Elec. J. Lin. Alg., 26 (2013)   673--691.

 \bibitem{AAB} J. Ahn, C. Alar, B. Bjorkman, S. Butler, J. Carlson, A. Goodnight, H. Knox, C. Monroe, M.C. Wigal, Ordered multiplicity inverse eigenvalue problem for graphs on six vertices, Elec. J. Lin. Alg., 37 (2021) 316--358.



\bibitem{BFH} W. Barrett, S. Fallat, H.T. Hall, L. Hogben, J.C.-H. Lin, B.L. Shader, Generalizations of the strong Arnold property and the minimum number of distinct eigenvalues of a graph, Electron. J. Combin., 24(2) (2017) 1--28.

 \bibitem{private}
 W. Barrett, S. Fallat, V. Furst, S. Nasserasr, B. Rooney, M. Tait. Private Communication. 2022.
 
 \bibitem{BHL} W. Barrett, H. van der Holst, and R. Loewy, Graphs whose minimal rank is two, Electron. J. Lin. Alg., 11 (2004) 258--280.


 \bibitem{Ber} D. S. Bernstein, {\it Matrix Mathematics}, Princeton University Press, New York, 2005.





  \bibitem{Cav} M. S. Cavers, Clique partitions and coverings of graphs, Masters thesis, Waterloo, Ontario, Canada, 2005.

  \bibitem{CGMJ} Z. Chen, M. Grimm, P. McMichael, and C.R. Johnson, Undirected graphs of Hermitian matrices that admit only two distinct eigenvalues, Linear Algebra Appl., 458 (2014) 403--428.

 \bibitem{ct} V. Consonni, R. Todeschini, New spectral indices for molecule description, MATCH Commun. Math. Comput. Chem., 60 (2008) 3--14.

 \bibitem{cv} D. Cvetkovi\'c, P. Rowlinson,  S. Simi\'c, {\it An introduction to the theory of graph spectra\/}, Cambridge University Press, Cambridge, 2012.

 \bibitem{AH5} K. C. Das, S. A. Mojallal, Relation between signless Laplacian energy, energy of graph and its line graph, Linear Algebra Appl., 493 (2016) 91--107.

 \bibitem{AH6} K. C. Das, S. A. Mojallal, I. Gutman, Relations between Degrees, Conjugate Degrees and Graph Energies, Linear Algebra Appl., 515 (2017) 24--37.

 \bibitem{Erd} P. Erd\"{o}s, A. W. Goodman, L. P\'{o}sa, The representation of a graph by set intersections, Can. J. Math., 18 (1966) 106--112.


  \bibitem{FH} S. Fallat, L. Hogben, The minimum rank of symmetric matrices described by a graph: A survey, Linear Algebra Appl., 426 (2007) 558–582.


  \bibitem{Fal} S. Fallat, L. Hogben, Variants on the minimum rank problem: A survey II, Preprint, arXiv:1102-5142v1, 2011.

  \bibitem{Fer}  W. E. Ferguson, The construction of Jacobi and periodic Jacobi matrices with prescribed spectra, Math. Comp., 35 (1980) 1203--1220.


 \bibitem{GR} C. Godsil, G. Royle, {\it Algebraic Graph Theory\/}, Springer, New York, 2001.

 \bibitem{GU1} I. Gutman, Bounds for total $\pi$-electron energy of polymethines, Chem. Phys. Lett., 50 (1977) 488--490.

  \bibitem{GU2} I. Gutman, Bounds for all graph energies, Chem. Phys. Lett., 528 (2012) 72--74.

  \bibitem{GRMC} I. Gutman, M. Robbiano, E. A. Martins, D. M. Cardoso, L. Medina, O. Rojo, Energy of line graphs, Linear Algebra Appl., 433 (2010) 1312--1323.


 \bibitem{Har}  F. Harary, {\em Graph Theory}, Addison-Wesley, Reading, 1972.


 \bibitem{Hog} L. Hogben, Spectral graph theory and the inverse eigenvalue problem of a graph, Electron. J. Lin. Alg., 14 (2005) 12--31.

 

  \bibitem{KM} J. Koolen, V. Moulton, Maximal energy graphs, Adv. Appl. Math., 26 (2001) 47--52.


 \bibitem{LOS} R. H. Levene, P. Oblak, H. \v{S}migoc, A Nordhaus–Gaddum conjecture for the minimum number of distinct eigenvalues of a graph,  Linear Algebra Appl., 564 (2019) 236--263.

\bibitem{LOS2} R. H. Levene, P. Oblak, H. \v{S}migoc, Orthogonal symmetric matrices and joins of graphs,
Linear Algebra Appl., 652 (2022), 213–238.
 

  \bibitem{LC} B.-J. Li, G. J. Chang, Clique coverings and partitions of line graphs, Discrete Math., 308 (2008) 2075--2079.

  \bibitem{LSG} X. Li, Y. Shi, I. Gutman, {\em Graph energy}, Springer, New York, 2012.

 \bibitem{MC} B. J. McClelland, Properties of the latent roots of a matrix: The estimation of $\pi$-electron energies. J. Chem. Phys., 54. (1971) 640--643.

 \bibitem{MCR} S. McGuinness, R. Rees,  On the number of distinct minimal clique partitions and clique covers of a line graph, Discrete Math., 83 (1990) 49--62.





 \bibitem{Ni-1}  V. Nikiforov, The energy of graphs and matrices, J.Math. Anal. Appl., 326 (2007) 1472--1475.

 \bibitem{Ni-2} V. Nikiforov, Graphs and matrices with maximal energy,  J. Math. Anal., Appl. 327 (2007) 735--738.

 \bibitem{Ni-3}  V. Nikiforov, Extremal norms of graphs and matrices,  J. Math. Sci., 182 (2012) 164--174.


 \bibitem{RT}
 K.B. Reid and C. Thomassen, Edge sets contained in circuits, Israel J. of Math., 24 (1976) 305--319.

 \bibitem{so} W. So, M. Robbiano, N. M. de Abreu, I. Gutman, Applications of a theorem by Ky Fan in the theory of graph energy, Linear Algebra Appl.,
 432 (2010) 2163--2169.


 \end{thebibliography}
 \end{document}